%
%
%
%
%
\newif\ifsect\newif\iffinal
\secttrue\finaltrue
\def\strutdepth{\dp\strutbox}

\def\lsimb#1{\vadjust{\vtop to0pt{\baselineskip\strutdepth\vss
	\llap{\ttt\string #1\ }\null}}}
\def\rsimb#1{\vadjust{\vtop to0pt{\baselineskip\strutdepth\vss
	\line{\kern\hsize\rlap{\ttt\ \string #1}}\null}}}
\def\ssect #1. {\bigbreak\indent{\bf #1.}\enspace\message{#1}}
\def\smallsect #1. #2\par{\bigbreak\noindent{\bf #1.}\enspace{\bf #2}\par
	\global\parano=#1\global\eqnumbo=1\global\thmno=1\global\esno=0
	\global\remno=0\global\remno=0\global\defno=0
	\nobreak\smallskip\nobreak\noindent\message{#2}}
\def\thm #1: #2{\medbreak\noindent{\bf #1:}\if(#2\thmp\else\thmn#2\fi}
\def\thmp #1) { (#1)\thmn{}}
\def\thmn#1#2\par{\enspace{\sl #1#2}\par
        \ifdim\lastskip<\medskipamount \removelastskip\penalty 55\medskip\fi}
\def\qedn{\thinspace\null\nobreak\hfill\hbox{\vbox{\kern-.2pt\hrule height.2pt 
depth.2pt\kern-.2pt\kern-.2pt \hbox to2.5mm{\kern-.2pt\vrule width.4pt
\kern-.2pt\raise2.5mm\vbox to.2pt{}\lower0pt\vtop to.2pt{}\hfil\kern-.2pt
\vrule width.4pt\kern-.2pt}\kern-.2pt\kern-.2pt\hrule height.2pt depth.2pt
\kern-.2pt}}\par\medbreak}
\def\pf{\ifdim\lastskip<\smallskipamount \removelastskip\smallskip\fi
        \noindent{\sl Proof\/}:\enspace}
\def\itm#1{\par\indent\llap{\rm #1\enspace}\ignorespaces}

\def\forclose#1{\hfil\llap{$#1$}\hfilneg}
\def\newforclose#1{
	\ifsect\xdef #1{(\number\parano.\number\eqnumbo)}\else
	\xdef #1{(\number\eqnumbo)}\fi
	\hfil\llap{$#1$}\hfilneg
	\global \advance \eqnumbo by 1
	\iffinal\else\rsimb#1\fi}
\def\newforclosea#1{
	\ifsect\xdef #1{{\rm(\the\Apptok.\number\eqnumbo)}}\else
	\xdef #1{(\number\eqnumbo)}\fi
	\hfil\llap{$#1$}\hfilneg
	\global \advance \eqnumbo by 1
	\iffinal\else\rsimb#1\fi}
\def\forevery#1#2$${\displaylines{\let\eqno=\forclose\let\neweqa=\newforclosea
        \let\neweq=\newforclose\hfilneg\rlap{$\qqquad\forall#1$}\hfil#2\cr}$$}
\newcount\parano
\newcount\eqnumbo
\newcount\thmno
\newcount\versiono
\versiono=0
\def\neweqt#1$${\xdef #1{(\number\parano.\number\eqnumbo)}
	\eqno #1$$
	\iffinal\else\rsimb#1\fi
	\global \advance \eqnumbo by 1}
\def\newthmt#1 #2: #3{\xdef #2{\number\parano.\number\thmno}
	\global \advance \thmno by 1
	\medbreak\noindent
	\iffinal\else\lsimb#2\fi
	{\bf #1 #2:}\if(#3\thmp\else\thmn#3\fi}
\def\neweqf#1$${\xdef #1{(\number\eqnumbo)}
	\eqno #1$$
	\iffinal\else\rlap{$\smash{\hbox{\hfilneg\string#1\hfilneg}}$}\fi
	\global \advance \eqnumbo by 1}
\def\newthmf#1 #2: #3{\xdef #2{\number\thmno}
	\global \advance \thmno by 1
	\medbreak\noindent
	\iffinal\else\llap{$\smash{\hbox{\hfilneg\string#1\hfilneg}}$}\fi
	{\bf #1 #2:}\if(#3\thmp\else\thmn#3\fi}
\def\inizia{\ifsect\let\neweq=\neweqt\else\let\neweq=\neweqf\fi
\ifsect\let\newthm=\newthmt\else\let\newthm=\newthmf\fi}
\def\bititolo{\empty}
\gdef\begin #1 #2\par{\xdef\titolo{#2}
\ifsect\let\neweq=\neweqt\else\let\neweq=\neweqf\fi
\ifsect\let\newthm=\newthmt\else\let\newthm=\newthmf\fi
\centerline{\titlefont\titolo}
\if\bititolo\empty\else\medskip\centerline{\titlefont\bititolo}
\xdef\titolo{\titolo\ \bititolo}\fi
\bigskip
\centerline{\bigfont \autore}
\if\istituto!\else\bigskip
\centerline{\istituto}
\centerline{\indirizzo}
\centerline{\email}\fi
\medskip
\centerline{#1~\anno}
\bigskip\bigskip
\ifsect\else\global\thmno=1\global\eqnumbo=1\fi}
\font\titlefont=cmssbx10 scaled \magstep1
\font\bigfont=cmr12

\font\ttt=cmtt10 at 10truept
\font\eightrm=cmr8

\let\sc=\smallcaps

\font\bbr=msbm10
\font\sbbr=msbm7 
\font\ssbbr=msbm5
\def\ca #1{{\cal #1}}
\nopagenumbers
\binoppenalty=10000
\relpenalty=10000
\newfam\amsfam
\textfont\amsfam=\bbr \scriptfont\amsfam=\sbbr \scriptscriptfont\amsfam=\ssbbr
\let\de=\partial
\def\eps{\varepsilon}
\def\phe{\varphi}

\def\End{\mathop{\rm End}\nolimits}

\def\Re{\mathop{\rm Re}\nolimits}
\def\Im{\mathop{\rm Im}\nolimits}

\def\id{\mathop{\rm id}\nolimits}
\mathchardef\void="083F
\def\Z{{\mathchoice{\hbox{\bbr Z}}{\hbox{\bbr Z}}{\hbox{\sbbr Z}}
{\hbox{\sbbr Z}}}}
\def\R{{\mathchoice{\hbox{\bbr R}}{\hbox{\bbr R}}{\hbox{\sbbr R}}
{\hbox{\sbbr R}}}}
\def\C{{\mathchoice{\hbox{\bbr C}}{\hbox{\bbr C}}{\hbox{\sbbr C}}
{\hbox{\sbbr C}}}}
\def\N{{\mathchoice{\hbox{\bbr N}}{\hbox{\bbr N}}{\hbox{\sbbr N}}
{\hbox{\sbbr N}}}}
\def\P{{\mathchoice{\hbox{\bbr P}}{\hbox{\bbr P}}{\hbox{\sbbr P}}
{\hbox{\sbbr P}}}}
\def\Q{{\mathchoice{\hbox{\bbr Q}}{\hbox{\bbr Q}}{\hbox{\sbbr Q}}
{\hbox{\sbbr Q}}}}

\def\qqquad{\quad\qquad}

\newcount\notitle
\notitle=1
\headline={\ifodd\pageno\rhead\else\lhead\fi}
\def\rhead{\ifnum\pageno=\notitle\iffinal\hfill\else\hfill\tt Version 
\the\versiono; \the\day/\the\month/\the\year\fi\else\hfill\eightrm\titolo\hfill
\folio\fi}
\def\lhead{\ifnum\pageno=\notitle\hfill\else\eightrm\folio\hfill\autore\hfill
\fi}
\def\autore{Marco Abate}
\output={\plainoutput}
\newbox\bibliobox
\def\setref #1{\setbox\bibliobox=\hbox{[#1]\enspace}
	\parindent=\wd\bibliobox}
\def\biblap#1{\noindent\hang\rlap{[#1]\enspace}\indent\ignorespaces}
\def\art#1 #2: #3! #4! #5 #6 #7-#8 \par{\biblap{#1}#2: {\sl #3\/}.
	#4 {\bf #5} (#6)\if.#7\else, \hbox{#7--#8}\fi.\par\smallskip}
\def\book#1 #2: #3! #4 \par{\biblap{#1}#2: {\bf #3.} #4.\par\smallskip}
\def\coll#1 #2: #3! #4! #5 \par{\biblap{#1}#2: {\sl #3\/}. In {\bf #4,} 
#5.\par\smallskip}
\def\pre#1 #2: #3! #4! #5 \par{\biblap{#1}#2: {\sl #3\/}. #4, #5.\par\smallskip}
%
%
\let\newthm=\newthmt
\let\neweq=\neweqt
\newcount\esno\newcount\defno\newcount\remno
\def\Def #1\par{\global \advance \defno by 1
    \medbreak
{\bf Definition \the\parano.\the\defno:}\enspace #1\par
\ifdim\lastskip<\medskipamount \removelastskip\penalty 55\medskip\fi}
\def\Rem #1\par{\global \advance \remno by 1
    \medbreak
{\bf Remark \the\parano.\the\remno:}\enspace #1\par
\ifdim\lastskip<\medskipamount \removelastskip\penalty 55\medskip\fi}
\def\Es #1\par{\global \advance \esno by 1
    \medbreak
{\sc Example \the\parano.\the\esno:}\enspace #1\par
\ifdim\lastskip<\medskipamount \removelastskip\penalty 55\medskip\fi}
\def\istituto{Dipartimento di Matematica, Universit\`a di Pisa}
\def\indirizzo{Via Buonarroti 2, 56127 pisa}
\def\anno{2003}
\def\email{E-mail: abate@dm.unipi.it}
%
%
%
%
\begin {June} Discrete local holomorphic dynamics

\smallsect 1. Introduction

Let $M$ be a complex manifold, and $p\in M$. In this survey, a {\sl
(discrete) holomorphic local dynamical system} at~$p$ will be a holomorphic map
$f\colon U\to M$ such that~$f(p)=p$, where $U\subseteq M$ is an open
neighbourhood of~$p$; we shall also assume that $f\not\equiv\id_U$. We shall
denote by~$\End(M,p)$ the set of holomorphic local dynamical systems at~$p$. 

\Rem Since we are mainly concerned with the behavior of~$f$ nearby~$p$, we shall
sometimes replace~$f$ by its restriction to some suitable open
neighbourhood of~$p$. It is possible to formalize this fact by using germs of
maps and germs of sets at~$p$, but for our purposes it will be enough
to use a somewhat less formal approach.

\Rem In this survey we shall never have the occasion of discussing continuous
holomorphic dynamical systems (i.e., holomorphic foliations). So from now on all
dynamical systems in this paper will be discrete, except where explicitely noted
otherwise.

To talk about the dynamics of an $f\in\End(M,p)$ we need to define the iterates
of~$f$. If $f$ is defined on the set~$U$, then the second iterate $f^2=f\circ f$
is defined on~$U\cap f^{-1}(U)$ only, which still is an open neighbourhood
of~$p$. More generally, the $k$-th iterate $f^k=f\circ f^{k-1}$ is defined
on~$U\cap f^{-1}(U)\cap\cdots\cap f^{-(k-1)}(U)$. Thus it is natural to
introduce the {\sl stable set}~$K_f$ of~$f$ by setting
$$
K_f=\bigcap_{k=0}^\infty f^{-k}(U).
$$
Clearly, $p\in K_f$, and so the stable set is never empty (but it can happen
that $K_f=\{p\}$; see the next section for an example). The stable set of~$f$ is
the set of all points~$z\in U$ such that the {\sl orbit} $\{f^k(z)\mid
k\in\N\}$ is well-defined. If $z\in U\setminus K_f$, we shall say that $z$ (or
its orbit) {\sl escapes} from~$U$. 

Thus the first natural question in local holomorphic dynamics is:
\smallskip
\item{(Q1)} {\it What is the topological structure of~$K_f$?}
\smallskip
\noindent For instance, when does $K_f$ have non-empty interior? As we shall see
in section~4, holomorphic local dynamical systems such that $p$ belongs to the
interior of the stable set enjoy special properties; we shall then say that $p$
is {\sl stable} for~$f\in\End(M,p)$ if it belongs to the interior of~$K_f$.

\Rem Both the definition of stable set and Question~1 (as well as several other
definitions or questions we shall meet later on) are topological in character;
we might state them for local dynamical systems which are continuous only. As
we shall see, however, the {\it answers} will strongly depend on the
holomorphicity of the dynamical system.

Clearly, the stable set~$K_f$ is {\sl completely $f$-invariant,} that is
$f^{-1}(K_f)=K_f$ (this implies, in particular, that $f(K_f)\subseteq K_f$).
Therefore the pair $(K_f,f)$ is a discrete dynamical system in the usual
sense, and so the second natural question in local holomorphic dynamics is
\smallskip
\item{(Q2)} {\it What is the dynamical structure of~$(K_f,f)$?}
\smallskip
\noindent For instance, what is the asymptotic behavior of the orbits? Do
they converge to~$p$, or have they a chaotic behavior? Is there a dense
orbit? Do there exist proper {\sl $f$-invariant} subsets, that is sets
$L\subset K_f$ such that~$f(L)\subseteq L$? If they do exist, what is the
dynamics on them?

To answer all these questions, the most efficient way is to replace $f$ by a
``dynamically equivalent" but simpler (e.g., linear) map~$g$. In our context,
``dynamically equivalent" means ``locally conjugated"; and we have at least
three kinds of conjugacy to consider.

Let $f_1\colon U_1\to M_1$ and $f_2\colon U_2\to M_2$ be two holomorphic local
dynamical systems at~$p_1\in M_1$ and~$p_2\in M_2$ respectively. We shall say
that~$f_1$ and~$f_2$ are {\sl holomorphically} (respectively, {\sl
topologically\/}) {\sl locally conjugated} if there are open neighbourhoods
$W_1\subseteq U_1$ of~$p_1$, $W_2\subseteq U_2$ of~$p_2$, and a
biholomorphism (respectively, a homeomorphism) $\phe\colon W_1\to W_2$ with
$\phe(p_1)=p_2$ such that 
$$
f_1=\phe^{-1}\circ f_2\circ\phe
\qquad\hbox{on}\qquad\phe^{-1}\bigl(
W_2\cap f_2^{-1}(W_2)\bigr)=W_1\cap f_1^{-1}(W_1).
$$  
In particular we have
$$
\forevery{k\in\N}\qquad\qquad f_1^k=\phe^{-1}\circ f_2^k\circ\phe
\quad\hbox{on}\quad\phe^{-1}\bigl(
W_2\cap\cdots\cap f_2^{-(k-1)}(W_2)\bigr)=W_1\cap\cdots\cap f_1^{-(k-1)}(W_1),
$$
and thus $K_{f_2|_{W_2}}=\phe(K_{f_1|_{W_1}})$. So the local dynamics of~$f_1$
about~$p_1$ is to all purposes equivalent to the local dynamics of~$f_2$
about~$p_2$.

\Rem Using local coordinates centered at~$p\in M$ it is easy to show that any
holomorphic local dynamical system at~$p$ is holomorphically locally conjugated
to a holomorphic local dynamical system at~$O\in\C^n$, where $n=\dim M$.

Whenever we have an equivalence relation in a class of objects, there are
obvious classification problems. So the third natural question in local
holomorphic dynamics is
\smallskip
\item{(Q3)} {\it Find a (possibly small) class $\ca F$ of holomorphic local
dynamical systems at~$O\in\C^n$ such that every holomorphic local dynamical
system~$f$ at a point in an $n$-dimensional complex manifold is holomorphically
(respectively, topologically) locally conjugated to a (possibly) unique element
of~$\ca F$, called the {\sl holomorphic} (respectively, {\sl topological\/})
{\sl normal form} of~$f$.}
\smallskip
\noindent Unfortunately, the holomorphic classification is often too complicated
to be practical; the family~$\ca F$ of normal forms might be uncountable. A
possible replacement is looking for invariants instead of normal forms:
\smallskip
\item{(Q4)} {\it Find a way to associate a (possibly small) class of (possibly
computable) objects to any holomorphic local dynamical system~$f$ at~$O\in\C^n$,
called the {\sl invariants} of~$f$, so that two holomorphic local dynamical
systems at~$O$ can be holomorphically conjugated only if they have the
same invariants. The class of invariants is furthermore said {\sl complete} if
two holomorphic local dynamical systems at~$O$ are holomorphically conjugated
if and only if they have the same invariants.}
\smallskip
\noindent As remarked before, up to now all the questions we asked make sense
for topological local dynamical systems; the next one instead makes sense only
for holomorphic local dynamical systems.

A holomorphic local dynamical system at~$O\in\C^n$ is clearly given by an
element of~$\C_0\{z_1,\ldots,z_n\}^n$, the space of $n$-uples of converging
power series in~$z_1,\ldots,z_n$ without constant terms. The
space~$\C_0\{z_1,\ldots,z_n\}^n$ is a subspace of the
space~$\C_0[[z_1,\ldots,z_n]]^n$ of $n$-uples of formal power series without
constant terms. An element $\Phi\in \C_0[[z_1,\ldots,z_n]]^n$ has an inverse
(with respect to composition) still belonging to~$\C_0[[z_1,\ldots,z_n]]^n$ if
and only if its linear part is a linear automorphism of~$\C^n$. We shall 
say that two holomorphic local dynamical systems
$f_1$,~$f_2\in\C_0\{z_1,\ldots,z_n\}^n$ are {\sl formally conjugated} if there
exists an invertible $\Phi\in\C_0[[z_1,\ldots,z_n]]^n$ such that
$f_1=\Phi^{-1}\circ f_2\circ\Phi$ in~$\C_0[[z_1,\ldots,z_n]]^n$. 

It is clear that two holomorphically locally conjugated holomorphic local
dynamical systems are both formally and topologically locally conjugated too. On
the other hand, we shall see examples of holomorphic local dynamical systems
that are topologically locally conjugated without being neither formally nor
holomorphically locally conjugated, and examples of holomorphic local
dynamical systems that are formally conjugated without being
neither holomorphically nor topologically locally conjugated. So the last
natural question in local holomorphic dynamics we shall deal with is
\smallskip
\item{(Q5)}{\it Find normal forms and invariants with respect to the relation
of formal conjugacy for holomorphic local dynamical systems at~$O\in\C^n$.}
\smallskip
\noindent In this survey we shall present some of the main results known on
these questions, starting from the one-dimensional situation. But before
entering the main core of this paper I would like to heartily thank Mohamad
Pouryayevali for the wonderful and very warm hospitality I had the pleasure
to enjoy during my stay in Iran.

\smallsect 2. One complex variable: the hyperbolic case

Let us then start by discussing holomorphic local dynamical systems at~$0\in\C$.
As remarked in the previous section, such a system
is given by a converging power series~$f$ without constant term:
$$
f(z)=a_1z+a_2z^2+a_3z^3+\cdots\in\C_0\{z\}.
$$
The number $a_1=f'(0)$ is the {\sl multiplier} of~$f$.

Since $a_1 z$ is the best linear approximation of~$f$, it is sensible to expect
that the local dynamics of~$f$ will be strongly influenced by the value of~$a_1$.
For this reason we introduce the following definitions:
\smallskip
\item{--} if $|a_1|<1$ we say that the fixed point $0$ is {\sl attracting;}
\item{--} if $a_1=0$ we say that the fixed point $0$ is {\sl superattracting;}
\item{--} if $|a_1|>1$ we say that the fixed point $0$ is {\sl repelling;}
\item{--} if $|a_1|\ne 0$,~1 we say that the fixed point $0$ is {\sl
hyperbolic;}
\item{--} if $a_1\in S^1$ is a root of unity, we say that the fixed point~$0$ is
{\sl parabolic} (or {\sl rationally indifferent\/});
\item{--} if $a_1\in S^1$ is not a root of unity, we say that the fixed 
point~$0$ is {\sl elliptic} (or {\sl irrationally indifferent\/}).
\smallskip
\noindent As we shall see in a minute, the dynamics of one-dimensional
holomorphic local dynamical systems with a hyperbolic fixed point is pretty
elementary; so we start with this case. Notice that if 0 is an attracting (we
shall discuss the superattracting case momentarily) fixed point
for~$f\in\End(\C,0)$, then it is a repelling fixed point for the inverse
map~$f^{-1}\in\End(\C,0)$. 

Assume first that $0$ is attracting for the holomorphic local
dynamical system~$f\in\End(\C,0)$. Then we can write
$f(z)=a_1z+O(z^2)$, with~$0<|a_1|<1$; hence we can find a large constant~$C>0$,
a small constant~$\eps>0$ and
$0<\delta<1$ such that if $|z|<\eps$ then
$$
|f(z)|\le (|a_1|+C\eps)|z|\le\delta|z|.
\neweq\eqduuno
$$
In particular, if~$\Delta_\eps$ denotes the disk of center~$0$ and
radius~$\eps$, we have $f(\Delta_\eps)\subset\Delta_\eps$ for $\eps>0$ small
enough, and the stable set of~$f|_{\Delta_\eps}$ 
is~$\Delta_\eps$ itself (in particular, an one-dimensional attracting fixed point
is always stable). Furthermore,
$$
|f^k(z)|\le\delta^k|z|\to 0
$$
as $k\to+\infty$, and thus every orbit starting in~$\Delta_\eps$ is attracted by
the origin, which is the reason of the name ``attracting" for such a fixed
point. 

If instead 0 is a repelling fixed point, a similar argument (or the observation
that 0 is attracting for~$f^{-1}$) shows that for~$\eps>0$ small enough the
stable set of~$f|_{\Delta_\eps}$ reduces to the origin only: all (non-trivial)
orbits escape.

It is also not difficult to find holomorphic and topological normal forms for
one-dimensional holomorphic local dynamical systems with a hyperbolic fixed
point, as shown in the
following result, which marked the beginning of the theory of holomorphic
dynamical systems:

\newthm Theorem \Koenigs: (K\oe nigs, 1884 [K\oe]) Let $f\in\End(\C,0)$ be an
one-dimensional holomorphic local dynamical system with a hyperbolic fixed
point at the origin, and let $a_1\in\C^*$ be its multiplier. Then:
{\smallskip
\itm{(i)} $f$ is
holomorphically (and hence formally) locally conjugated to its linear
part~$g(z)=a_1 z$. 
\item{\rm (ii)} Two such holomorphic local dynamical systems are
holomorphically conjugated if and only if they have the same multiplier.
\item{\rm (iii)} $f$ is topologically locally conjugated to the
map $g_<(z)=z/2$ if $|a_1|<1$, and to the map $g_>(z)=2z$\break
\indent if~$|a_1|>1$.}

\ifdim\lastskip<\smallskipamount \removelastskip\smallskip\fi
\noindent{\sl Sketch of proof\/}:\enspace Let us assume $0<|a_1|<1$; if
$|a_1|>1$ it will suffice to apply the same argument to~$f^{-1}$.

Put $\phe_k=f^k/a_1^k$; using \eqduuno\ it is not difficult to show that the
sequence $\{\phe_k\}$ converges to a holomorphic map
$\phe\colon\Delta_\eps\to\C$ for $\eps>0$ small enough. Since $\phe_k'(0)=1$ for
all~$k\in\N$, we have $\phe'(0)=1$ and so, up to possibly shrink~$\eps$, we can
assume that $\phe$ is a biholomorphism with its image. Moreover, we have
$$
\phe\bigl(f(z)\bigr)=\lim_{k\to+\infty}{f^k\bigl(f(z)\bigr)\over a_1^k}
=a_1\lim_{k\to+\infty}{f^{k+1}(z)\over a_1^{k+1}}=a_1\phe(z),
$$
that is $f=\phe^{-1}\circ g\circ \phe$, as claimed.

Since $f_1=\phe^{-1}\circ f_2\circ\phe$ implies $f_1'(0)=f_2'(0)$, the
multiplier is invariant under holomorphic local conjugation, and so two
one-dimensional holomorphic local dynamical systems with a hyperbolic fixed point
are holomorphically locally conjugated if and only if they have the same
multiplier.

Finally, if $|a_1|<1$ it is easy to build a topological conjugacy between $g$
and~$g_<$ on~$\Delta_\eps$: it suffices to choose any homeomorphism~$\phe$
between the annulus $\{\eps/2\le |z|<\eps\}$ and the annulus $\{|a_1|\eps\le
|z|<\eps\}$, and to extend it by induction to a homeomorphism between the annuli
$\{\eps/2^k\le |z|\le\eps/2^{k-1}\}$ and $\{|a_1|^k\eps\le
|z|\le|a_1|^{k-1}\eps\}$ by requiring 
$$
\phe({\textstyle{1\over2}}z)=a_1\,\phe(z).
$$
Putting finally $\phe(0)=0$ we then get the topological conjugacy we were
looking for.
\qedn

Notice that $g_<(z)={1\over 2}z$ and $g_>(z)=2z$ cannot be topologically
conjugated, because (for instance) the origin is stable for~$g_<$ and it is not
stable for~$g_>$.

Thus the dynamics in the one-dimensional hyperbolic case is completely clear.
The superattracting case can be treated similarly. If $0$ is a superattracting
point for an~$f\in\End(\C,0)$, we can write 
$$
f(z)=a_rz^r+a_{r+1}z^{r+1}+\cdots
$$
with $a_r\ne 0$; the number $r\ge 2$ is the {\sl order} of the superattracting
point. An argument similar to the one described above shows that for $\eps>0$
small enough the stable set of~$f|_{\Delta_\eps}$ still is all of~$\Delta_\eps$,
and the orbits converge (faster than in the attracting case) to the origin.
Furthermore, replacing the maps~$\phe_k$ in the proof of Theorem~\Koenigs\ by
maps of the form
$$
\phe_k(z)=[f^k(z)]^{1/r^k},
$$
for a suitable choice of the $r^k$-th root, one can prove the following

\newthm Theorem \Bottcher: (B\"ottcher, 1904 [B]) Let $f\in\End(\C,0)$ be an
one-dimensional holomorphic local dynamical system with a superattracting fixed
point at the origin, and let $r\ge 2$ be its order. Then:
{\smallskip
\itm{(i)} $f$ is
holomorphically (and hence formally) locally conjugated to the map~$g(z)=z^r$.
\itm{(ii)} two such holomorphic local dynamical systems are holomorphically
(or topologically) conjugated if and\break\indent only if they have the same
order.}

Therefore the one-dimensional local dynamics about a hyperbolic or
superattracting fixed point is completely clear; let us now discuss what happens
about a parabolic fixed point.

\smallsect 3. One complex variable: the parabolic case

Let $f\in\End(\C,0)$ be a (non-linear) holomorphic local dynamical system with a
parabolic fixed point at the origin. Then we can write 
$$
f(z)=e^{2i\pi p/q}z+a_{r+1} z^{r+1}+a_{r+2}z^{r+2}+\cdots,
\neweq\eqtuno
$$
with $a_{r+1}\ne 0$, where $p/q\in\Q\cap[0,1)$ is the {\sl rotation number}
of~$f$, and the number $r+1\ge 2$ is the {\sl multiplicity} of~$f$ at the fixed
point. 

The first observation is that such a dynamical system is never locally
conjugated to its linear part, not even topologically, unless it is of finite
order. Indeed, if we had
$\phe^{-1}\circ f\circ\phe(z)=e^{2\pi ip/q}z$ we would have $\phe^{-1}\circ
f^q\circ\phe=\id$, that is~$f^q=\id$. 

In particular, if the rotation number is~0 (that is the multiplier is~1, and
we shall say that $f$ is {\sl tangent to the identity\/}), then $f$ {\it
cannot} be locally conjugated to the identity (unless it was the identity
to begin with, which is not a very interesting case dynamically
speaking). More precisely, the stable set of such an
$f$ is never a neighbourhood of the origin. To understand why, let us first
consider a map of the form
$$
f(z)=z(1+az^r)
$$
for some $a\ne 0$. Let $v\in S^1\subset\C$ be such that $av^r$ is real and
positive. Then for any $c>0$ we have
$$
f(cv)=c(1+c^rav^r)v\in\R^+v;
$$
moreover, $|f(cv)|>|cv|$. In other words, the half-line~$\R^+v$ is $f$-invariant
and repelled from the origin, that is $K_f\cap\R^+ v=\void$. Conversely, if
$av^r$ is real and negative then the segment $[0,|a|^{-1/r}]v$ is $f$-invariant
and attracted by the origin. So $K_f$ neither is a neighbourhood of the origin
nor reduces to~$\{0\}$. 

This example suggests the following definition. Let $f\in\End(\C,0)$ be of the
form~\eqtuno\ and tangent to the identity. Then a unit vector~$v\in S^1$
is an {\sl attracting} (respectively, {\sl repelling\/}) {\sl direction}
for~$f$ at the origin if~$a_{r+1}v^r$ is real and negative (respectively,
positive). Clearly, there are $r$ equally spaced attracting directions,
separated by $r$ equally spaced repelling directions; furthermore, a repelling
(attracting) direction for~$f$ is attracting (repelling) for~$f^{-1}$, which is
defined in a neighbourhood of the origin.

It turns out that to every attracting direction is associated a connected
component of~$K_f\setminus\{0\}$. Let~$v\in S^1$ be an attracting direction for
an~$f$ tangent to the identity. The {\sl basin} centered at~$v$ is the set of
points~$z\in K_f\setminus\{0\}$ such that $f^k(z)\to 0$ and $f^k(z)/|f^k(z)|\to
v$ (notice that, up to shrinking the domain of~$f$, we can assume that $f(z)\ne
0$ for all $z\in K_f\setminus\{0\}$). If $z$ belongs to the basin centered
at~$v$, we shall say that the orbit of~$z$ {\sl tends to~$0$ tangent to~$v$.}

A slightly more specialized (but more useful) object is the following: an {\sl
attracting petal} centered at an attracting direction~$v$ is an open simply
connected $f$-invariant set $P\subseteq K_f\setminus\{0\}$ such that a point
$z\in K_f\setminus\{0\}$ belongs to the basin centered at~$v$ if and only if its
orbit intersects~$P$. In other words, the orbit of a point tends to~$0$ tangent
to~$v$ if and only if it is eventually contained in~$P$. A {\sl repelling petal}
(centered at a repelling direction) is an attracting petal for the inverse
of~$f$.

It turns out that the basins centered at the attracting directions are exactly
the connected components of~$K_f\setminus\{0\}$, as shown in the {\it Leau-Fatou
flower theorem:} 

\newthm Theorem \flower: (Leau, 1897 [L]; Fatou, 1919-20 [F1--3]) Let
$f\in\End(\C,0)$ be a holomorphic local dynamical system tangent to the identity
with multiplicity~$r+1\ge 2$ at the fixed point. Let $v_1,v_3,\ldots,v_{2r-1}\in
S^1$ be the
$r$ attracting directions of~$f$ at the origin, and $v_2,v_4,\ldots,v_{2r}\in
S^1$ the $r$ repelling directions. Then {\smallskip
\item{\rm (i)} There exists for each attracting (repelling) direction~$v_{2j-1}$
($v_{2j}$) an attracting (repelling) petal~$P_{2j-1}$ ($P_{2j}$), so that the
union of these $2r$ petals together with the origin forms a neighbourhood
of the origin. Furthermore, the
$2r$ petals are arranged ciclically so that two petals intersect if and only if
the angle between their central directions is~$\pi/r$.
\item{\rm (ii)} $K_f\setminus\{0\}$ is the (disjoint) union of the basins
centered at the
$r$ attracting directions.
\item{\rm (iii)} If $P$ is an attracting petal, then $f|_P$ is holomorphically
conjugated to the translation $z\mapsto z+1$ defined\break\indent on a subset of
the complex plane containing some right half-plane.}

\ifdim\lastskip<\smallskipamount \removelastskip\smallskip\fi
\noindent{\sl Sketch of proof\/}:\enspace Up to a linear change of variables, we
can assume that $a_{r+1}=-1$, so that the attracting directions are the $r$-th
roots of unity. For any $\delta>0$, the set $\{z\in\C\mid|z^r-\delta|<\delta\}$
has exactly $r$ connected components, each one centered on a different $r$-th
root of unity; it will turns out that, for $\delta$ small enough, these connected
components are the attracting petals of~$f$.

Let $P_\delta$ denote one of these connected components, and let $\psi\colon
P_\delta\to\C$ be given by
$$
\psi(z)={1\over rz^r}.
$$
This is a biholomorphism of~$P_\delta$ with a right half-plane
$H_\delta=\{w\in\C\mid\Re w>1/(2r\delta)\}$, and we have
$$
\psi\circ f\circ\psi^{-1}(w)=w+1+O(w^{-1/r}).
\neweq\eqtdue
$$
Then, setting $F=\psi\circ f\circ\psi^{-1}$, it is not difficult to prove that
for $\delta>0$ small enough the right half-plane~$H_\delta$ is $F$-invariant,
and that for any $w\in H_\delta$ the orbit $\{F^k(w)\}$ converges to~$\infty$
tangent to~$+1$. Thus it follows that~$P_\delta$ is $f$-invariant, and that the
orbits in~$P_\delta$ tends to the origin tangent to the central direction~$v$
of~$P_\delta$. Since every orbit converging to the origin tangent to~$v$ must
eventually intersect~$P_\delta$, every such~$P_\delta$ is an attracting petal.

Arguing in the same way with~$f^{-1}$ we get the repelling petals, and thus (i)
follows. Since it is not difficult to prove that every orbit converging to the
origin must be tangent to an attracting direction, (ii) follows too.
Finally, a subtler argument shows that we can modify~$\psi$ in each petal so to
get rid of the term~$O(w^{-1/r})$ in~\eqtdue, proving (iii).\qedn

So we have a complete description of the dynamics in the neighbourhood of the
origin. Actually, Camacho has pushed this argument even further, obtaining a
complete topological classification of one-dimensional holomorphic local
dynamical systems tangent to the identity:

\newthm Theorem \Camacho: (Camacho, 1978 [C]; Shcherbakov, 1982 [S]) Let
$f\in\End(\C,0)$ be a holomorphic local dynamical system tangent to the identity
with multiplicity~$r+1$ at the fixed point. Then $f$ is topologically locally
conjugated to the map
$$
z\mapsto z+z^{r+1}.
$$

The formal classification is simple too, though different, and it can be
obtained with an easy computation (see, e.g., Milnor [Mi]):

\newthm Proposition \formaltangent: Let $f\in\End(\C,0)$ be a
holomorphic local dynamical system tangent to the identity with
multiplicity~$r+1$ at the fixed point. Then $f$ is formally conjugated to the map
$$
z\mapsto z+z^{r+1}+\beta z^{2r+1},
$$
where $\beta$ is a formal (and holomorphic) invariant given
by
$$
\beta={1\over 2\pi i}\int_\gamma{dz\over z-f(z)},
\neweq\eqttre
$$
where the integral is taken over a small positive loop~$\gamma$ about the
origin.

The number~$\beta$ given by~\eqttre\ is called {\sl index} of~$f$ at the fixed
point.

The holomorphic classification is much more complicated: as shown by
Voronin~[V] and \'Ecalle [\'E1--2] in~1981, it depends on functional invariants. 
We shall now try to
roughly describe it; see~[I2] (and the original
papers; see also~[K]) for details. Let $f\in\End(\C,0)$ be tangent to the
identity with multiplicity~$r+1$ at the fixed point; up to a linear
change of coordinates we can assume that~$a_{r+1}=1$. Let
$P_1,\ldots,P_{2r}$ be a set of petals as in Theorem~\flower.(i), chosen so
that $P_{2r}$ is centered on the positive real semiaxis, and the others are
arranged cyclically counterclockwise. Denote by
$H_j$ the biholomorphism conjugating $f|_{P_j}$ to the shift~$z\mapsto z+1$ in
either a right (if $j$ is odd) or left (if $j$ is even) half-plane given by
Theorem~\flower.(iii) --- applied to~$f^{-1}$ for the repelling petals. If we
moreover require that 
$$
H_j(z)=-{1\over rz^r}+\beta\log z+o(1),
\neweq\eqtquattro
$$
where $\beta$ is the index of~$f$ at the origin, then $H_j$ is uniquely
determined. Thus in the sets $H_j(P_j\cap P_{j+1})$ we can consider the
composition
$\tilde\Phi_j=H_{j+1}\circ H_j^{-1}$. It is easy to check that
$\tilde\Phi_j(w+1)=\tilde\Phi_j(w)+1$ for~$j=1,\ldots,2r-1$, and thus
$\psi_j=\tilde\Phi_j-\id$ is a 1-periodic holomorphic function (for
$j=2r$ we need to take~$\psi_{2r}=\Phi_{2r}=\id+2\pi i\beta$ to get a
1-periodic function). Hence each $\psi_j$ can be extended to a
suitable upper (if
$j$ is odd) or lower (if~$j$ is even) half-plane. Furthermore, it is possible to
prove that the functions~$\psi_1,\ldots,\psi_{2r}$ are exponentially
decreasing, that is they are bounded by $\exp(-c|w|)$ as $|\Im w|\to+\infty$, for
a suitable~$c>0$ depending on~$f$.

Now, if we replace $f$ by a holomorphic local conjugate~$g=h^{-1}\circ f\circ
h$, and denote by $G_j$ the corresponding biholomorphisms, it turns out that
$H_j\circ G_j^{-1}=\id+a$ for a suitable $a\in\C$ independent of~$j$. This
suggests the introduction of an equivalence relation on the set of $2r$-uple of
functions of the kind $(\psi_1,\ldots,\psi_{2r})$.

Let $M_r$ denote the set of $2r$-uple of holomorphic 1-periodic
functions~$\psi=(\psi_1,\ldots,\psi_{2r})$, with $\psi_j$ defined in a suitable
upper (if $j$ is odd) or lower (if $j$ is even) half-plane, and exponentially
decreasing when~$|\Im w|\to+\infty$. We shall say that $\psi$,~$\tilde\psi\in
M_r$ are {\sl equivalent} if there is~$a\in\C$ such that
$\tilde\psi_j=\psi_j\circ(\id+a)$ for $j=1,\ldots,2r$. We denote by~$\ca M_r$ the
set of all equivalence classes.

The procedure described above allows us to associate to any $f\in\End(\C,0)$
tangent to the identity with multiplicity~$r+1$ at the fixed point an
element~$\mu_f\in\ca M_r$, called the {\sl sectorial invariant.} Then the
holomorphic classification proved by
\'Ecalle and Voronin is

\newthm Theorem \EV: (\'Ecalle, 1981 [\'E1--2]; Voronin, 1981 [V]) Let
$f$,~$g\in\End(\C,0)$ be two holomorphic local dynamical systems tangent to the
identity. Then $f$ and $g$ are holomorphically locally conjugated if and only
if they have the same multiplicity, the same index and the same sectorial
invariant. Furthermore, for any $r\ge1$, $\beta\in\C$ and $\mu\in\ca M_r$ there
exists $f\in\End(\C,0)$ tangent to the identity with multiplicity~$r+1$,
index~$\beta$ and sectorial invariant~$\mu$.  

For a sketch of the proof, together with a more geometrical description of the
sectorial invariant, see~[I2] and [M1--2]. 

\Rem In particular, holomorphic local dynamical systems tangent to the
identity give examples of local dynamical systems that are topologically
conjugated without being neither holomorphically nor formally conjugated, and of
local dynamical systems that are formally conjugated without being
holomorphically conjugated.

Finally, if $f\in\End(\C,0)$ satisfies $a_1=e^{2\pi i p/q}$, then $f^q$ is
tangent to the identity. Therefore we can apply the previous results to~$f^q$
and then infer informations about the dynamics of the original~$f$. See~[Mi],
[C], [\'E1--2] and~[V] for details.

\smallsect 4. One complex variable: the elliptic case

We are left with the elliptic case:
$$
f(z)=e^{2\pi i\theta}z+a_2z^2+\cdots\in\C_0\{z\},
\neweq\eqquno
$$
with $\theta\notin\Q$. It turns out that the local dynamics depends mostly
on the numerical properties of~$\theta$. More precisely, for a full measure
subset~$B$ of $\theta\in[0,1]\setminus\Q$ all holomorphic local dynamical systems
of the form~\eqquno\ are {\sl holomorphically linearizable,} that is
holomorphically locally conjugated to their (common) linear part, the irrational
rotation
$z\mapsto e^{2\pi i\theta}z$. Conversely, the complement $[0,1]\setminus B$ is a
$G_\delta$-dense set, and for all~$\theta\in [0,1]\setminus B$ the quadratic
polynomial~$z\mapsto z^2+e^{2\pi i\theta}z$ is not holomorphically linearizable.
This is the gist of the results due to Cremer, Siegel, Bryuno and Yoccoz we are
going to describe in this section.

The first worthwhile observation in this setting is that it is possible to give
a topological characterization of the holomorphically linearizable local
dynamical systems:

\newthm Proposition \local: Let $f\in\End(\C,0)$ be a holomorphic local
dynamical system with multiplier $0<|\lambda|\le 1$. Then
$f$ is holomorphically linearizable if and only if it is topologically
linearizable if and only if $0$ is stable for~$f$.

\ifdim\lastskip<\smallskipamount \removelastskip\smallskip\fi
\noindent{\sl Sketch of proof\/}:\enspace Assume that $0$ is stable. If
$0<|\lambda|<1$, we already saw that $f$ is linearizable. If $|\lambda|=1$, set
$$
\phe_k(z)={1\over k}\sum_{j=0}^{k-1}{f^j(z)\over\lambda^j},
$$
so that
$$
\phe_k\circ f=\lambda\phe_{k+1}+{\lambda\over k}(\phe_{k+1}-f).
\neweq\eqqdue
$$
The stability of~$0$ implies that $\{\phe_k\}$ is a normal family in a
neighbourhood of the origin, and \eqqdue\ implies that a converging subsequence
converges to a conjugation between~$f$ and the rotation~$z\mapsto\lambda z$.\qedn

The second important observation is that two elliptic holomorphic local dynamical
systems with the same multiplier are always formally conjugated:

\newthm Proposition \qdue: Let $f\in\End(\C,0)$ be a holomorphic local dynamical
system of multiplier $\lambda=e^{2\pi i\theta}\in S^1$ with $\theta\notin\Q$.
Then $f$ is formally conjugated to its linear part.

\ifdim\lastskip<\smallskipamount \removelastskip\smallskip\fi
\noindent{\sl Sketch of proof\/}:\enspace It is an easy computation to prove
that there is a unique formal power series 
$$
h(z)=z+h_2z^2+\cdots\in\C[[z]]
$$ 
such that
$h(\lambda z)=f\bigl(h(z)\bigr)$. For later use we explicitely remark that the
coefficients of the formal linearization satisfy
$$
h_j={a_j+X_j\over\lambda^j-\lambda},
\neweq\eqqtre
$$
where $X_j$ is a polynomial expression in
$a_2,\ldots,a_{j-1},h_2,\ldots,h_{j-1}$.\qedn

The formal power series linearizing~$f$ is not converging if its coefficients
grow too fast. Thus \eqqtre\ links the radius of convergence of~$h$ to the
behavior of~$\lambda^j-\lambda$: if the latter becomes too small, the series
defining~$h$ does not converge. This is known as the {\it small denominators
problem} in this context. 

It is then natural to introduce the following quantity:
$$
\Omega_\lambda(m)=\min_{1\le k\le m}|\lambda^k-1|,
$$
for $\lambda\in S^1$ and $m\ge 1$. Clearly, $\lambda$ is a root of unity if and
only if $\Omega_\lambda(m)=0$ for all $m$ greater or equal to some~$m_0\ge 1$;
furthermore, 
$$
\lim_{m\to+\infty}\Omega_\lambda(m)=0
$$ 
for all~$\lambda\in S^1$. 

The first one to actually prove that there are elliptic holomorphic local
dynamical systems not linearizable has been Cremer, in 1927 [Cr1]. Later he
proved the following:

\newthm Theorem \Cremer: (Cremer, 1938 [Cr2]) Let $\lambda\in S^1$ be
such that
$$
\limsup_{m\to+\infty}\left(-{1\over m}\log\Omega_\lambda(m)\right)=+\infty.
\neweq\eqqquattro
$$
Then there exists $f\in\End(\C,0)$ with multiplier~$\lambda$ which is not
holomorphically linearizable.
Furthermore, the set of $\lambda\in S^1$ satisfying~$\eqqquattro$ contains
a~$G_\delta$-dense set. 

\ifdim\lastskip<\smallskipamount \removelastskip\smallskip\fi
\noindent{\sl Sketch of proof\/}:\enspace Choose inductively $a_j\in\{0,1\}$ so
that $|a_j+X_j|\ge1/2$ for all~$j\ge 2$, where~$X_j$ is as in~\eqqtre. Then
$$
f(z)=\lambda z+a_2z^2+\cdots\in\C_0\{z\}
$$
while \eqqquattro\ implies that the radius of convergence of the formal
linearization~$h$ is~0, and thus $f$ cannot be holomorphically linearizable, as
required.

Finally, let $S(q_0)\subset S^1$ denote the set of $\lambda=e^{2\pi i\theta}\in
S^1$ such that
$$
\left|\theta-{p\over q}\right|<{1\over 2^{q!}}
$$
for some~$p/q\in\Q$ in lowest terms with $q\ge q_0$. Then it is not difficult to
check that each $S(q_0)$ is a dense open set in~$S^1$, and that all
$\lambda\in\bigcap_{q_0\ge 1}S(q_0)$ satisfy~\eqqquattro.\qedn

On the other hand, Siegel, using the technique of majorant series, in 1942 gave
a condition on the multiplier ensuring holomorphic linearizability:

\newthm Theorem \Siegel: (Siegel, 1942 [Si]) Let $\lambda\in S^1$ be such
that there exists $\beta\ge 1$ and $\gamma>0$ such that
$$
\forevery{m\ge 2}{1\over\Omega_\lambda(m)}\le\gamma\, m^\beta.
\neweq\eqqcinque
$$
Then all $f\in\End(\C,0)$ with multiplier~$\lambda$ are
holomorphically linearizable.
Furthermore, the set of $\lambda\in S^1$ satisfying~$\eqqcinque$ for
some~$\beta\ge 1$ and $\gamma>0$ is of full Lebesgue measure in~$S^1$.

\Rem It is interesting to notice that for generic (in a topological sense)
$\lambda\in S^1$ there is a non-linearizable holomorphic local dynamical system
with multiplier~$\lambda$, while for almost all (in a measure-theoretic sense)
$\lambda\in S^1$ every holomorphic local dynamical system
with multiplier~$\lambda$ is holomorphically linearizable.

A bit of terminology is now useful: if $f\in\End(\C,0)$ is elliptic, we shall
say that the origin is a {\sl Siegel point} if $f$ is holomorphically 
linearizable; otherwise it is a {\sl Cremer point.}

Theorem~\Siegel\ suggests the existence of a number-theoretical condition
on~$\lambda$ ensuring that the origin is a Siegel point for any holomorphic
local dynamical system of multiplier~$\lambda$. And indeed this is the content
of the celebrated {\sl Bryuno-Yoccoz theorem:}

\newthm Theorem \BY: Let
$\lambda\in S^1$. {\smallskip 
\item{\rm(i)} {\rm (Bryuno, 1965 [Bry1--3])} If $\lambda$
satisfies 
$$
\sum_{k=0}^{+\infty}\left(- 2^{-k}\log\Omega_\lambda(2^{k+1})\right)<+\infty,
\neweq\eqqsei
$$
then the origin is a Siegel point for all $f\in\End(\C,0)$ with
multiplier~$\lambda$.
\item{\rm (ii)} {\rm (Yoccoz, 1988 [Y1--2])} If $\lambda$ does not satisfy
$\eqqsei$, then the origin is a Cremer point for some 
$f\in\End(\C,0)$\break\indent with multiplier~$\lambda$.
In particular, the origin is a Cremer point for~$f(z)=\lambda z+z^2$.}

The original proof by Bryuno of Theorem~\BY.(i) uses majorant series; see, e.g.,
[He] and references therein. Yoccoz found a more geometric approach, based on
conformal and quasi-conformal geometry, and
proved~Theorem~\BY.(ii). Furthermore, he showed that the origin is a Siegel
point for all elliptic holomorphic local dynamical systems with
multiplier~$\lambda$ if and only if it is a Siegel point for~$f(z)=\lambda
z+z^2$. See also [P9].

\Rem Condition \eqqsei\ is usually expressed in a different way. Write
$\lambda=e^{2\pi i\theta}$, and let $\{p_k/q_k\}$ be the sequence of rational
numbers converging to~$\theta$ given by the expansion in continued fractions.
Then \eqqsei\ is equivalent to
$$
\sum_{k=0}^{+\infty} {1\over q_k}\log q_{k+1}<+\infty,
$$
while \eqqcinque\ is equivalent to~$q_{n+1}=O(q_n^\beta)$, and \eqqquattro\
is equivalent to
$$
\limsup_{k\to+\infty} {1\over q_k}\log q_{k+1}=+\infty.
$$
See [He], [Y2] and references therein for details.

If $0$ is a Siegel point for~$f\in\End(\C,0)$, the local dynamics of~$f$ is
completely clear, and simple enough. On the other hand, if 0 is a Cremer point
of~$f$, then the local dynamics of $f$ is very complicated and not yet
completely understood. P\'erez-Marco (in [P2, 4--7]) has studied the topology and
the dynamics of the stable set in this case. Some of his results are summarized
in the following

\newthm Theorem \PerezMarco: (P\'erez-Marco, 1995 [P6, 7]) Assume that $0$ is
a Cremer point for an elliptic holomorphic local dynamical system
$f\in\End(\C,0)$. Then:
{\smallskip
\item{\rm(i)} The stable set $K_f$ is compact, connected, full (i.e.,
$\C\setminus K_f$ is connected), it is not reduced to~$\{0\}$, and it is not
locally connected at any point distinct from the origin.
\item{\rm(ii)} Any point of~$K_f\setminus\{0\}$ is recurrent (that is, 
a limit point of its orbit).
\item{\rm(iii)} There is an orbit in~$K_f$ which accumulates at the origin, but
no non-trivial orbit converges to the origin.}

\Rem As far as I know, there are neither a topological nor a holomorphic
complete classification of elliptic holomorphic dynamical systems with a Cremer
point. Furthermore, if $\lambda\in S^1$ is not a root of unity and does not
satisfy Bryuno's condition~\eqqsei, we can find $f_1$,~$f_2\in\End(\C,0)$ with
multiplier~$\lambda$ such that~$f_1$ is holomorphically linearizable while~$f_2$
is not. Then $f_1$ and $f_2$ are formally conjugated without being neither
holomorphically nor topologically locally conjugated.

See also [P1, 3] for other results on the dynamics about a Cremer point.

\smallsect 5. Several complex variables: the hyperbolic case

Now we start the discussion of local dynamics in several complex variables. In
this case the theory is much less complete than its one-variable counterpart.

Let $f\in\End(\C^n,O)$ be a holomorphic local dynamical system at $O\in\C^n$,
with $n\ge 2$. We can write $f$ using a {\sl homogeneous expansion}
$$
f(z)=P_1(z)+P_2(z)+\cdots\in\C_0\{z_1,\ldots,z_n\}^n,
$$
where $P_j$ is an $n$-uple of homogeneous polynomials of degree~$j$. In
particular,~$P_1$ is the differential~$df_O$ of~$f$ at the origin, and $f$ is
locally invertible if and only if~$P_1$ is invertible.

We have seen that in dimension one the multiplier (i.e., the derivative at the
origin) plays a main r\^ole. When $n>1$, a similar r\^ole is played by the
eigenvalues of the differential. Thus we introduce the following definitions:
\smallskip
\item{--} if all eigenvalues of $df_O$ have modulus less than 1, we say that the
fixed point~$O$ is {\sl attracting;}
\item{--} if all eigenvalues of $df_O$ have modulus greater than 1, we say that
the fixed point~$O$ is {\sl repelling;}
\item{--} if all eigenvalues of $df_O$ have modulus different from 1, we say that
the fixed point~$O$ is {\sl hyperbolic} (notice that we allow the eigenvalue
zero);
\item{--} if all eigenvalues of $df_O$ are roots of unity, we say that the
fixed point~$O$ is {\sl parabolic;} in particular, if $df_O=\id$ we say that $f$
is {\sl tangent to the identity;}
\item{--} if all eigenvalues of $df_O$ have modulus 1 but none is a root of
unity, we say that the fixed point~$O$ is {\sl elliptic;}
\item{--} if $df_O=O$, we say that the fixed point~$O$ is {\sl superattracting.}
\smallskip
\noindent Other cases are clearly possible, but for our aims this list is
enough. In this survey we shall be mainly concerned with hyperbolic and
parabolic fixed points; however, in the last section we shall also present some
results valid in other cases.

Let us begin assuming that the origin is a hyperbolic
fixed point for an $f\in\End(\C^n,O)$ not necessarily invertible. We then have a
canonical splitting
$$
\C^n=E^s\oplus E^u,
$$
where $E^s$ (respectively, $E^u$) is the direct sum of the generalized
eigenspaces associated to the eigenvalues of~$df_O$ with modulus less
(respectively, greater) than 1. Then the first main result in this subject is
the famous {\sl stable manifold theorem} (originally due to Perron~[Pe] and
Hadamard~[H]; see [FHY, HK, HPS, Pes, Sh] for proofs in
the $C^\infty$ category, Wu~[Wu] for a proof in the holomorphic category, and
[A3] for a proof in the non-invertible case):

\newthm Theorem \stable: Let $f\in\End(\C^n,O)$ be a holomorphic local dynamical
system with a hyperbolic fixed point at the origin. Then:
{\smallskip
\item{\rm(i)}the stable set~$K_f$ is an embedded complex submanifold of (a
neighbourhood of the origin in)~$\C^n$, tangent to~$E^s$ at the origin;
\item{\rm(ii)}there is an embedded complex submanifold~$W_f$ of (a
neighbourhood of the origin in)~$\C^n$, called the\break\indent {\it unstable
set} of~$f$, tangent to~$E^u$ at the origin, such that $f|_{W_f}$ is
invertible, $f^{-1}(W_f)\subseteq W_f$, and $z\in W_f$\break\indent if and only
if there is a sequence $\{z_{-k}\}_{k\in\N}$ in the domain of~$f$ such that
$z_0=z$ and $f(z_{-k})=z_{-k+1}$ for\break\indent all~$k\ge 1$. Furthermore, if
$f$ is invertible then $W_f$ is the stable set of~$f^{-1}$. }

The proof is too involved to be summarized here; it suffices to say that both
$K_f$ and $W_f$ can be recovered, for instance, as fixed points of a suitable
contracting operator in an infinite dimensional space (see the references quoted
above for details).

\Rem If the origin is an attracting fixed point, then $E^s=\C^n$, and $K_f$ is an
open neighbourhood of the origin, its {\sl basin of attraction.} However, as we
shall discuss below, this does not imply that
$f$ is holomorphically linearizable, not even when it is invertible.
Conversely, if the origin is a repelling fixed point, then
$E^u=\C^n$, and~$K_f=\{O\}$. Again, not all holomorphic local dynamical systems
with a repelling fixed point are holomorphically linearizable.

If a point in the domain~$U$ of a holomorphic local dynamical system with a
hyperbolic fixed point does not belong either to the stable set or to the
unstable set, it escapes both in forward time (that is, its orbit escapes) and in
backward time (that is, it is not the end point of an infinite orbit contained
in~$U$). In some sense, we can think of the stable and unstable sets (or, as
they are usually called in this setting, stable and unstable {\it manifolds\/})
as skewed coordinate planes at the origin, and the orbits outside these
coordinate planes follow some sort of hyperbolic path, entering and leaving any
neighbourhod of the origin in finite time.

Actually, this idea of straightening stable and unstable manifolds can be
brought to fruition (at least in the invertible case), and it yields one of the
possible proofs (see~[HK, Sh, A3] and references therein) of the
{\sl Grobman-Hartman theorem:}

\newthm Theorem \GrobmanHartman: (Grobman, 1959 [G1--2]; Hartman, 1960 [Har])
Let
$f\in\End(\C^n,O)$ be a locally invertible holomorphic local dynamical system
with a hyperbolic fixed point. Then $f$ is topologically locally conjugated to
its differential~$df_O$.

Thus, at least from a topological point of view, the local dynamics about an
invertible hyperbolic fixed point is completely clear. This is definitely not the
case if the local dynamical system is not invertible in a neighbourhood of the
fixed point. For instance, already Hubbard and Papadopol [HP] noticed that
a B\"ottcher-type theorem for superattracting points in several complex
variables is just not true: there are holomorphic local dynamical systems with a
superattracting fixed point which are not even topologically locally conjugated
to the first non-vanishing term of their homogeneous expansion. Recently, Favre
and Jonsson~[FJ] have begun a very detailed study of superattracting fixed
points in~$\C^2$, study that should lead to their topological classification.

The holomorphic and even the formal classification are not as simple as the
topological one. The main problem is that, if we denote by
$\lambda_1,\ldots,\lambda_n\in\C$ the eigenvalues of~$df_O$, then it may happen
that
$$
\lambda_1^{k_1}\cdots\lambda_n^{k_n}-\lambda_j=0
\neweq\eqcuno
$$
for some $1\le j\le n$ and some $k_1,\ldots,k_n\in\N$ with $k_1+\cdots+k_n\ge
2$; a relation of this kind is called a {\sl resonance} of~$f$. When $n=1$ there
is a resonance if and only if the multiplier is a root of unity, or
zero; but if
$n>1$ resonances may occur in the hyperbolic case too. Anyway, a computation
completely analogous to the one yielding Proposition~\qdue\ proves the following

\newthm Proposition \qtre: Let $f\in\End(\C^n,O)$ be a (locally invertible)
holomorphic local dynamical system with a hyperbolic fixed point and no
resonances. Then $f$ is formally conjugated to its differential~$df_O$.

In presence of resonances, even the formal classification is not that easy.
Let us assume, for simplicity, that $df_O$ is in Jordan form, that is
$$
P_1(z)=(\lambda_1z,\epsilon_2 z_1+\lambda_2z_2,\ldots,\epsilon_n
z_{n-1}+\lambda_nz_n),
$$
with $\epsilon_1,\ldots,\epsilon_{n-1}\in\{0,1\}$. We shall say that a
monomial
$z_1^{k_1}\cdots z_n^{k_n}$ in the $j$-th coordinate of~$f$ is {\sl resonant} if
$k_1+\cdots+k_n\ge 2$ and
$\lambda_1^{k_1}\cdots\lambda_n^{k_n}=\lambda_j$. Then the Proposition~\qtre\
can be generalized to

\newthm Proposition \PoincareDulac: Let $f\in\End(\C^n,O)$ be a locally 
invertible holomorphic local dynamical system with a hyperbolic fixed point.
Then it is formally conjugated to a $g\in\C_0[[z_1,\ldots,z_n]]^n$ such that
$dg_O$ is in Jordan normal form, and $g$ has only resonant monomials.

The formal series~$g$ is called {\sl Poincar\'e-Dulac normal form} of~$f$; see
Arnold~[Ar] for a proof of Proposition~\PoincareDulac. 

The problem with Poincar\'e-Dulac normal forms is that they are not unique. In
particular, one may wonder whether it could be possible to have such a normal
form including {\it finitely many} resonant monomials only (as happened, for
instance, in Proposition~\formaltangent). This is
indeed the case (see, e.g., Reich~[R1]) when $df_O$ belongs to the so-called {\sl
Poincar\'e domain,} that is when $df_O$ is invertible and $O$ is either
attracting or repelling (when
$df_O$ is still invertible but does not belong to the Poincar\'e domain, we
shall say that it belongs to the {\sl Siegel domain\/}). As far as I know, the
problem of finding canonical formal normal forms 
when
$df_O$ belongs to the Siegel domain (and $f$ is hyperbolic) is still open.

It should be remarked that, in the hyperbolic case, the problem of formal
linearization is equivalent to the problem of smooth linearization. This has
been proved by Sternberg~[St1--2] and Chaperon~[Ch]:

\newthm Theorem \Sternberg: (Sternberg, 1957 [St1--2]; Chaperon, 1986 [Ch])
Let 
$f$,~$g\in\End(\C^n,O)$ be two holomorphic local dynamical
systems, and assume that $f$ is locally invertible and with a hyperbolic fixed
point at the origin. Then $f$ and $g$ are formally conjugated if and only if
they are smoothly locally conjugated. In particular, $f$ is smoothly
linearizable if and only if it is formally linearizable. Thus if there are
no resonances then $f$ is smoothly linearizable.

Even without resonances, the holomorphic linearizability is not guaranteed.
The easiest positive result is due to Poincar\'e~[Po] who, using majorant
series, proved the following

\newthm Theorem \Poincare: (Poincar\'e, 1893 [Po])  Let $f\in\End(\C^n,O)$ be
a  locally invertible holomorphic local dynamical system with an attracting or
repelling fixed point. Then $f$ is holomorphically linearizable if and
only if it is formally linearizable. In particular, if there are no resonances
then $f$ is holomorphically linearizable.

Reich~[R2] describes holomorphic normal forms when $df_O$ belongs to the
Poincar\'e domain and there are resonances (see also~[\'EV]); P\'erez-Marco [P8]
discusses the problem of holomorphic linearization in the presence of resonances.

When $df_O$ belongs to the Siegel domain, even without resonances, the formal
linearization might diverge. To describe the known results, let us introduce the
following quantity:
$$
\Omega_{\lambda_1,\ldots,\lambda_n}(m)=\min\bigl\{|\lambda_1^{k_1}\cdots
\lambda_n^{k_n}-\lambda_j|\bigm| k_1,\ldots,k_n\in\N,\,2\le k_1+\cdots+k_n\le m, 
\,1\le j\le n\bigr\}
$$
for $m\ge 2$ and $\lambda_1,\ldots,\lambda_n\in\C$. In particular, if
$\lambda_1,\ldots,\lambda_n$ are the eigenvalues of~$df_O$, we shall
write~$\Omega_f(m)$ for~$\Omega_{\lambda_1,\ldots,\lambda_n}(m)$.

It is clear that $\Omega_f(m)\ne 0$ for all~$m\ge 2$ if and only if there are no
resonances. It is also not difficult to prove that if $df_O$ belongs to the
Siegel domain then
$$
\lim_{m\to+\infty}\Omega_f(m)=0,
$$
which is the reason why, even without resonances, the formal linearization might
be diverging, exactly as in the one-dimensional case. As far as I know, the best
positive and negative results in this setting are due to Bryuno~[Bry2--3], and
are a natural generalization of their one-dimensional counterparts:

\newthm Theorem \Bryunohyp: (Bryuno, 1971 [Bry2--3]) Let $f\in\End(\C^n,O)$ be a
holomorphic local dynamical system such that
$df_O$ belongs to the Siegel domain, is linearizable and has no resonances.
Assume moreover that 
$$
\sum_{k=0}^{+\infty}\left(-{1\over 2^k}\log\Omega_f(2^{k+1})\right)<+\infty.
\neweq\eqcdue
$$
Then $f$ is holomorphically linearizable.

\newthm Theorem \BryunoCremer: Let $\lambda_1,\ldots,\lambda_n\in\C$ be without
resonances and such that
$$
\limsup_{m\to+\infty}\left(-{1\over
m}\log\Omega_{\lambda_1,\ldots,\lambda_n}(m)\right)=+\infty.
$$
Then there exists $f\in\End(\C^n,O)$, with $df_O=\hbox{\rm
diag}(\lambda_1,\ldots,\lambda_n)$, not holomorphically 
linearizable.

\Rem These theorems hold even without hyperbolicity assumptions. 

It should be remarked that, contrarily to the one-dimensional case, it is not
known whether condition \eqcdue\ is necessary for the holomorphic 
linearizability of all holomorphic local dynamical systems with a given linear
part belonging to the Siegel domain. See also P\"oschel~[P\"o] for a
generalization of~Theorem \Bryunohyp, and Il'yachenko~[I1] for an important
result related to Theorem~\BryunoCremer. Finally, in [DG] are discussed results
in the spirit of Theorem~\Bryunohyp\ without assuming that the differential is
diagonalizable.

\smallsect 6. Several complex variables: the parabolic case

A first natural question in the several complex variables parabolic case is
whether a result like the Leau-Fatou flower theorem holds, and, if so, in
which form. To present what is known on this subject in this section we shall
restrict our attention to holomorphic local dynamical systems 
tangent to the identity; consequences on dynamical systems with a more general
parabolic fixed point can be deduced taking a suitable iterate (but see also
the end of this section for results valid when the differential at the
fixed point is not diagonalizable).

So we are interested in the local dynamics of a holomorphic local dynamical
system $f\in\End(\C^n,O)$ of the form
$$
f(z)=z+P_\nu(z)+P_{\nu+1}(z)+\cdots\in\C_0\{z_1,\ldots,z_n\}^n,
$$
where $P_\nu$ is the first non-zero term in the homogeneous expansion of~$f$;
the number~$\nu\ge 2$ is the {\sl order} of~$f$. 

The two main ingredients in the statement of the Leau-Fatou flower theorem were
the attracting directions and the petals. Let us first describe
a several variables analogue of attracting directions.

Let $f\in\End(\C^n,O)$ be tangent at the identity and of order~$\nu$. A {\sl
characteristic direction} for~$f$ is a non-zero vector $v\in\C^n\setminus\{O\}$
such that $P_\nu(v)=\lambda v$ for some~$\lambda\in\C$. If $P_\nu(v)=O$ (that
is, $\lambda=0$) we shall say that $v$ is a {\sl degenerate} characteristic
direction; otherwise, (that is, if $\lambda\ne 0$) we shall say that $v$ is {\sl
non-degenerate.}

There is an equivalent definition of characteristic directions that shall be
useful later on. The $n$-uple of $\nu$-homogeneous polynomial~$P_\nu$ induces a
meromorphic self-map of~$\P^{n-1}(\C)$, still denoted by~$P_\nu$. Then, under
the canonical projection $\C^n\setminus\{O\}\to\P^{n-1}(\C)$ that we shall
denote by~$v\mapsto[v]$, the non-degenerate characteristic directions correspond
exactly to fixed points of~$P_\nu$, and the degenerate characteristic directions
correspond exactly to indeterminacy points of~$P_\nu$. By the way, using
Bezout's theorem it is easy to prove (see, e.g., [AT]) that the number of
characteristic directions of~$f$, counted according to a suitable multiplicity,
is given by~$(\nu^n-1)/(\nu-1)$. 

\Rem The characteristic directions are {\it complex} directions; in particular,
it is easy to check that $f$ and $f^{-1}$ have the same characteristic
directions. Later we shall see how to associate to (most) characteristic
directions $\nu-1$ petals, each one in some sense centered about a {\it
real} attracting direction corresponding to the same complex characteristic
direction.

The notion of characteristic directions has a dynamical origin. We shall say
that an orbit
$\{f^k(z_0)\}$ converges to the origin {\sl tangentially} to a
direction~$[v]\in\P^{n-1}(\C)$ if $f^k(z_0)\to O$ in~$\C^n$ and
$[f^k(z_0)]\to[v]$ in~$\P^{n-1}(\C)$. Then

\newthm Proposition \suno: Let $f\in\End(\C^n,O)$ be a holomorphic dynamical
system tangent to the identity. If there is an orbit of~$f$ converging to the
origin tangentially to a direction~$[v]\in\P^{n-1}(\C)$, then $v$ is a
characteristic direction of~$f$.

\ifdim\lastskip<\smallskipamount \removelastskip\smallskip\fi
\noindent{\sl Sketch of proof\/}:\enspace ([Ha2]) For simplicity let us
assume~$\nu=2$; a similar argument works for~$\nu>2$.  

If $v$ is a degenerate characteristic direction, there is nothing to prove. If
not, up to a linear change of coordinates we can write
$$
\cases{f_1(z)=z_1+p^1_2(z_1,z')+p_3^1(z_1,z')+\cdots,\cr
\noalign{\smallskip}
f'(z)=z'+p'_2(z_1,z')+p'_3(z_1,z')+\cdots,\cr}
$$
where $z'=(z_2,\ldots,z_n)\in\C^{n-1}$, $f=(f_1,f')$, $P_j=(p_j^1,p'_j)$ and so
on, with $v=(1,v')$ and $p_2^1(1,v')\ne 0$. Making the substitution 
$$
\cases{w_1=z_1,\cr
z'=w'z_1,\cr}
\neweq\eqblowup
$$
which is a change of variable outside the hyperplane $z_1=0$, the map~$f$ becomes
$$
\cases{\tilde f_1(w)=w_1+p^1_2(1,w')w_1^2+p_3^1(1,w')w_1^3+\cdots,\cr
\noalign{\smallskip}
\tilde f'(w)=w'+r(w')w_1+O(w_1^2),\cr}
\neweq\eqsdue
$$
where $r(w')$ is a polynomial such that $r(v')=O$ if and only if $(1,v')$ is a
characteristic direction of~$f$ with~$p_2^1(1,v')\ne0$.

Now, the hypothesis is that there exists an orbit $\{f^k(z_0)\}$ converging to
the origin and such that $[f^k(z_0)]\to[v]$. Writing $\tilde
f^k(w_0)=\bigl(w_1^k,(w')^k\bigr)$, this implies that $w_1^k\to 0$ and
$(w')^k\to v'$. Then it is not difficult to prove that 
$$
\lim_{k\to+\infty}{1\over kw_1^k}=-p_2^1(1,v')
$$
and then that $(w')^{k+1}-(w')^k$ is of the order of $r(v')/k$, which implies
$r(v')=O$, as claimed.\qedn

\Rem There are (unfortunately?) examples of $f\in\End(\C^2,O)$ tangent to the
identity with an orbit converging to the origin which is not tangent to any
direction (see [Ri1]).

The several variables analogue of a petal is instead given by the
notion of parabolic curve. A {\sl parabolic curve} for $f\in\End(\C^n,O)$
tangent to the identity is an injective holomorphic map
$\phe\colon\Delta\to\C^n\setminus\{O\}$ satisfying the following properties:
\smallskip
\item{(a)} $\Delta$ is a simply connected domain in~$\C$ with $0\in\de\Delta$;
\item{(b)} $\phe$ is continuous at the origin, and $\phe(0)=O$;
\item{(c)} $\phe(\Delta)$ is $f$-invariant, and $(f|_{\phe(\Delta)})^k\to O$
uniformly on compact subsets as $k\to+\infty$.
\smallskip
\noindent Furthermore, if $[\phe(\zeta)]\to[v]$ in $\P^{n-1}(\C)$ as~$\zeta\to
0$ in~$\Delta$, we shall say that the parabolic curve $\phe$ is {\sl tangent} to
the direction~$[v]\in\P^{n-1}(\C)$. 

Then the first main generalization of the Leau-Fatou flower theorem to several
complex variables is

\newthm Theorem \EcalleHakim: (\'Ecalle, 1985 [\'E3]; Hakim, 1998 [Ha2]) Let
$f\in\End(\C^n,O)$ be a holomorphic local dynamical system tangent to the
identity of order~$\nu\ge 2$. Then for any non-degenerate characteristic
direction~$[v]\in\P^{n-1}(\C)$ there exist (at least) $\nu-1$ parabolic curves
for~$f$ tangent to~$[v]$.

\ifdim\lastskip<\smallskipamount \removelastskip\smallskip\fi
\noindent{\sl Sketch of proof\/}:\enspace \'Ecalle proof is based on his theory
of resurgence of divergent series; we shall describe here the ideas behind
Hakim's proof, which depends on more standard arguments. 

For the sake of simplicity, let us assume $n=2$; without loss of generality we
can also assume $[v]=[1:0]$. Then after a linear change of variables and
a transformation of the kind~\eqblowup\ we are reduced
to prove the existence of a parabolic curve at the origin for a map of the form
$$
\cases{f_1(z)=z_1-z_1^\nu+O(z_1^{\nu+1},z_1^\nu z_2),\cr
\noalign{\smallskip}
f_2(z)=z_2\bigl(1-\lambda
z_1^{\nu-1}+O(z_1^\nu,z_1^{\nu-1} z_2)\bigr)+z_1^\nu
	\psi(z),\cr}
$$
where $\psi$ is holomorphic with $\psi(O)=0$, and $\lambda\in\C$.
Given
$\delta>0$, set
$D_{\delta,\nu}=\{\zeta\in\C\mid |\zeta^{\nu-1}-\delta|<\delta\}$; this open set
has $\nu-1$ connected components, all of them satisfying condition (a) on the
domain of a parabolic curve. Furthermore, if $u$ is a holomorphic function
defined on one of these connected components, of the form $u(\zeta)=\zeta^2
u_o(\zeta)$ for some bounded holomorphic function~$u_o$, and such that
$$
u\bigl(f_1\bigl(\zeta,u(\zeta)\bigr)\bigr)=f_2\bigl(\zeta,u(\zeta)\bigr),
\neweq\eqstre
$$
then it is not difficult to verify that $\phe(\zeta)=\bigl(\zeta,u(\zeta)\bigl)$
is a parabolic curve for~$f$ tangent to~$[v]$. 

So we are reduced to finding a solution of \eqstre\ in each connected component
of~$D_{\delta,\nu}$, with $\delta$ small enough. For any holomorphic $u=\zeta^2
u_o$ defined in such a connected component, let
$f_u(\zeta)=f_1\bigl(\zeta,u(\zeta)\bigr)$, put
$$
H(z)=z_2-{z_1^\lambda\over f_1(z)^\lambda} f_2(z),
$$
and define the operator $T$ by setting
$$
(Tu)(\zeta)=\zeta^\lambda \sum_{k=0}^\infty {H\bigl(f_u^k(\zeta), u\bigl(
f_u^k(\zeta)\bigr)\bigr)\over f_u^k(\zeta)^\lambda}.
$$
Then, if $\delta>0$ is small enough, it is possible to prove that $T$ is
well-defined, that $u$ is a fixed point of~$T$ if and only if it satisfies
\eqstre, and that $T$ is a contraction of a closed convex set of a
suitable complex Banach space --- and thus it has a fixed point. In this way if
$\delta>0$ is small enough we get a unique solution of \eqstre\ for each
connected component of~$D_{\delta,\nu}$, and hence $\nu-1$ parabolic curves
tangent to~$[v]$.\qedn

A set of $\nu-1$ parabolic curves obtained in this way will be called a {\sl
Fatou flower} for~$f$ tangent to~$[v]$. 

\Rem It should be remarked that a similar result for 2-dimensional maps with
$\lambda\notin\N^*$ has been obtained by Weickert [W] too; the
computations needed in the proof for the case $\lambda\in\N^*$ are considerably
harder, and were not carried out by him.

\Rem When there is a one-dimensional $f$-invariant complex submanifold passing
through the origin tangent to a characteristic direction~$[v]$, the
previous theorem is just a consequence of the usual one-dimensional theory. But
it turns out that in most cases such an $f$-invariant complex submanifold does
not exist: see [Ha2] for a concrete example, and [\'E3] for a general
discussion. 

We can also have $f$-invariant complex submanifolds of
dimension strictly greater than one still attracted by the origin. Given
a holomorphic local dynamical system $f\in\End(\C^n,O)$ tangent to the
identity of order~$\nu\ge 2$, and
a non-degenerate characteristic direction $[v]\in\P^{n-1}(\C)$, 
the eigenvalues $\alpha_1,\ldots,\alpha_{n-1}\in\C$ of the linear operator
$d(P_\nu)_{[v]}-\id\colon T_{[v]}\P^{n-1}(\C)\to T_{[v]}\P^{n-1}(\C)$
will be called the {\sl directors} of~$[v]$. Then, using a more elaborate
version of her proof of Theorem~\EcalleHakim, Hakim has been able to prove the
following:

\newthm Theorem \Hakim: (Hakim, 1997 [Ha3]) Let $f\in\End(\C^n,O)$ be a
holomorphic local dynamical system tangent to the identity of order~$\nu\ge 2$.
Let
$[v]\in\P^{n-1}(\C)$ be a non-degenerate characteristic direction, with directors
$\alpha_1,\ldots,\alpha_{n-1}\in\C$.
Furthermore, assume that $\Re\alpha_1,\ldots,\Re\alpha_d>0$ and
$\Re\alpha_{d+1},\ldots,\Re\alpha_{n-1}\le 0$ for a suitable $d\ge 0$. Then:
\smallskip
{\item{\rm(i)}There exists an $f$-invariant $(d+1)$-dimensional complex
submanifold~$M$ of~$\C^n$, with the origin in its boundary, such that the orbit
of every point of~$M$ converges to the origin tangentially to~$[v]$;
\item{\rm(ii)} $f|_M$ is holomorphically conjugated to the translation
$\tau(w_0,w_1,\ldots,w_d)=(w_0+1,w_1,\ldots,w_d)$ defined\break\indent on a
suitable right half-space in~$\C^{d+1}$. }

\Rem In particular, if all the directors of~$[v]$
have positive real part, there is an open domain attracted by the origin.
However, the condition given by Theorem~\Hakim\ is not necessary for the
existence of such an open domain; see Rivi~[Ri1] for an easy example, and
Ushiki~[Us] for a more elaborate example with an open domain attracted by the
origin where $f$ cannot be conjugate to a translation.

In his monumental work [\'E3] \'Ecalle has given a complete set of formal
invariants for holomorphic local dynamical systems tangent to the identity with
at least one non-degenerate characteristic direction. For instance, he
has proved the following

\newthm Theorem \Ecalleformal: (\'Ecalle, 1985 [\'E3]) Let $f\in\End(\C^n,O)$
be a holomorphic local dynamical system tangent to the identity of order~$\nu\ge
2$. Assume that
{\smallskip
\item{\rm(a)} $f$ has exactly $(\nu^n-1)/(\nu-1)$ distinct non-degenerate
characteristic directions and no degenerate characteristic directions;
\item{\rm(b)} the directors of any non-degenerate characteristic
direction are irrational and mutually independent over~$\Z$.
\smallskip
\noindent Choose a non-degenerate characteristic direction~$[v]\in\P^{n-1}(\C)$,
and let $\alpha_1,\ldots,\alpha_{n-1}\in\C$ be its directors. Then there exist a
unique $\rho\in\C$ and unique (up to dilations) formal
series~$R_1,\ldots,R_n\in\C[[z_1,\ldots,z_n]]$, where each $R_j$ contains only
monomial of total degree at least~$\nu+1$ and of partial degree in~$z_j$ at
most~$\nu-2$, such that $f$ is formally conjugated to the time-1 map of the
formal vector field
$$
X={1\over(\nu-1)(1+\rho z_n^{\nu-1})}\left\{[-z_n^\nu+R_n(z)]{\de\over\de
z_n}+\sum_{j=1}^{n-1}[-\alpha_j z_n^{\nu-1}z_j+R_j(z)]{\de\over\de z_j}\right\}.
$$}

Another approach to the formal classification, at least in dimension~2, is
described in~[BM].

Furthermore, using his theory of resurgence, and always assuming the
existence of at least one non-degenerate characteristic direction, \'Ecalle has
also provided a set of holomorphic invariants for holomorphic local dynamical
systems tangent to the identity, in terms of differential operators with formal
power series as coefficients. Moreover, if the directors of all
non-degenerate characteristic direction are irrational and satisfy a suitable
diophantine condition, then these invariants become a complete set of
invariants. See [\'E4] for a description of his results, and [\'E3]
for the details.

Now, all these results beg the question: what happens when there are no
non-degenerate characteristic directions? For instance, this is the case for
$$
\cases{f_1(z)=z_1+bz_1z_2+z_2^2,\cr
\noalign{\smallskip}
f_2(z)=z_2-b^2 z_1z_2-b z_2^2+z_1^3,\cr}
$$
for any $b\in\C^*$, and it is easy to build similar examples of any
order. At present, the theory in this case is satisfactorily developed for $n=2$
only. In particular, in [A2] is proved the following

\newthm Theorem \Abate: (Abate, 2001 [A2]) Every holomorphic local dynamical
system $f\in\End(\C^2,O)$ tangent to the identity, with an isolated fixed point,
admits at least one Fatou flower tangent to some direction.

\Rem Bracci and Suwa have proved a version of Theorem~\Abate\ for
$f\in\End(M,p)$ where $M$ is a {\it singular} variety with not too bad a
singularity at~$p$; see [BrS] for details.

Let us describe the main ideas in the proof of Theorem~\Abate, because they
provide some insight on the dynamical structure of holomorphic local
dynamical systems tangent to the identity, and on how to deal with it.

The first idea is to exploit in a systematic way the transformation~\eqblowup,
following a procedure borrowed from algebraic geometry. If $p$ is a point in a
complex manifold~$M$, there is a canonical way to build a complex
manifold~$\tilde M$, called the {\sl blow-up} of~$M$ at~$p$, provided with a
holomorphic projection $\pi\colon\tilde M\to M$, and such that $E=\pi^{-1}(p)$,
the {\sl exceptional divisor} of the blow-up, is canonically biholomorphic
to~$\P(T_pM)$, and $\pi|_{\tilde M\setminus E}\colon\tilde M\setminus E\to
M\setminus\{p\}$ is a biholomorphism. In suitable local coordinates, the
map~$\pi$ is exactly given by~\eqblowup. Furthermore, if $f\in\End(M,p)$ is
tangent to the identity, there is a unique way to lift~$f$ to a map~$\tilde
f\in\End(\tilde M,E)$ such that $\pi\circ\tilde f=f\circ\pi$, where $\End(\tilde
M,E)$ is the set of holomorphic maps defined in a neighbourhood of~$E$ with
values in~$\tilde M$ and which are the identity on~$E$. In particular, the
characteristic directions of~$f$ become points in the domain of~$\tilde f$.

This approach allows to determine which characteristic directions are
dynamically meaningful. Take $f=(f_1,f_2)\in\End(\C^2,O)$ tangent to the
identity; if $\ell=\gcd(f_1-z_1,f_2-z_2)$, we can write
$$
f_j(z)=z_j+\ell(z)g_j(z)
$$ 
for $j=1$,~2, with $g_1$ and $g_2$ relatively prime in~$\C\{z_1,z_2\}$. We shall
say that $O$ is a {\sl singular point} for~$f$ if~$g_1(O)=g_2(O)=0$. Clearly, if
$O$ is an isolated fixed point of~$f$ then it is singular; but if $O$ is not an
isolated fixed point (i.e., $\ell\not\equiv 1$) it might not be singular. Only
singular points are dynamically meaningful, because a not too difficult
computation (see~[A2], and [AT] for an $n$-dimensional generalization)
yields the following

\newthm Proposition \sing: Let $f\in\End(\C^2,O)$ be a holomorphic local
dynamical system tangent to the identity. If the fixed point~$O$ is not singular,
then $K_f$ reduces to the fixed point set of~$f$.

Now, if $\tilde M$ is the blow-up of~$\C^2$ at the origin, then the lift~$\tilde
f$ of $f$ belongs to~$\End(\tilde M,[v])$ for any direction~$[v]\in\P^1(\C)=E$.
We shall then say that $[v]\in\P^1(\C)$ is a {\sl singular direction} of~$f$ if
it is a singular point for~$\tilde f$. It turns out that non-degenerate
characteristic directions are always singular (but the converse does not
necessarily hold), and that singular directions are always characteristic (but
the converse does not necessarily hold): the singular directions are the
dynamically interesting characteristic directions.

The important feature of the blow-up procedure is that even though the underlying
manifold becomes more complex, the lifted maps become simpler. Indeed, using
an argument similar to one (described, for instance, in~[MM]) used in the study
of singular holomorphic foliations of 2-dimensional complex manifolds,
in~[A2] it is shown that after a finite number of blow-ups our original
holomorphic local dynamical system~$f\in\End(\C^2,O)$ can be lifted to a
map~$\tilde f$ whose singular points (are finitely many and) satisfy one of
the following conditions:
\smallskip
\item{(o)} they are {\sl dicritical,} that is with infinitely many singular
directions; or,
\item{($\star$)}in suitable local coordinates centered at the singular point we
can write
$$
\cases{\tilde f_1(z)=z_1+\ell(z)\bigl(\lambda_1 z_1+O(\|z\|^2)\bigr),\cr
\tilde f_2(z)=z_2+\ell(z)\bigl(\lambda_2 z_2+O(\|z\|^2)\bigr),
\cr}
$$
with
\itemitem{($\star_1$)} $\lambda_1$,~$\lambda_2\ne 0$ and
$\lambda_1/\lambda_2$,~$\lambda_2/\lambda_1\notin\N$, or
\itemitem{($\star_2$)} $\lambda_1\ne 0$, $\lambda_2=0$.

\Rem This ``reduction of the singularities" statement holds only in dimension~2,
and it is not clear how to replace it in higher dimensions.

It is not too difficult to prove that Theorem~\EcalleHakim\ (actually,
the ``easy" case of this theorem) can be applied both to dicritical and to
$(\star_1)$ singularities; therefore as soon as this blow-up procedure produces
such a singularity, we get a Fatou flower for the original dynamical system~$f$.

So to end the proof of Theorem~\Abate\ it suffices to prove that any such
blow-up procedure {\it must} produce at least one dicritical or $(\star_1)$
singularity. To get such a result, we need a completely new ingredient.

Let $E$ be a compact Riemann surface inside a 2-dimensional complex manifold
(for instance, $E$ can be the exceptional divisor of the blow-up of a point~$p$),
and take $f\in\End(M,E)$ tangent to the identity to all points of~$E$
(this happens, for instance, if $f$ is the lifting of a map tangent at the
identity at~$p$). Given~$q\in E$, choose local coordinates~$(z_1,z_2)$ in~$M$
centered at~$q$ and such that $E$ is locally given by~$\{z_2=0\}$. Then the
function
$$
k(z_1)=\lim_{z_2\to 0}{f_2(z)-z_2\over z_2(f_1(z)-z_1)}
$$
is either a meromorphic function defined in a neighbourhood of~$q$, or
identically~$\infty$. It turns out that:
\smallskip
\item{--} if $k$ is identically~$\infty$ at one point~$q\in E$, it is
identically~$\infty$ at all points of~$E$; in this case we shall say that $f$ is
{\sl not tangential} to~$E$;
\item{--} if $f$ is tangential to~$E$ (this happens, for instance, if $f$ is
obtained blowing up a non-dicritical singularity), then the residue of~$k$
at~$q$ is independent of the local coordinates used to define~$k$, and it is
called the {\sl index}~$\iota_q(f,E)$ of~$f$ at~$q$ along~$E$;
\item{--} if $f$ is tangential to~$E$, and $q\in E$ is not singular for~$f$, then
$\iota_q(f,E)=0$; in particular, $\iota_q(f,E)\ne 0$ only for a finite number of
points of~$E$.
\smallskip
\noindent Then following an argument suggested by Camacho and Sad [CS] in
their study of the separatrices of holomorphic foliations it is possible to
prove the following {\sl index theorem:}

\newthm Theorem \index: (Abate, 2001 [A2]) Let $E$ be a compact Riemann surface
inside a 2-dimensional complex manifold~$M$. Take $f\in\End(M,E)$ such that $f$
is tangent to the identity at all points of~$E$, and assume that~$f$ is
tangential to~$E$. Then
$$
\sum_{q\in E}\iota_q(f,E)=c_1(N_E),
$$
where $c_1(N_E)$ is the first Chern class of the normal bundle~$N_E$ of~$E$
in~$M$.

\Rem If $f$ is the lift to the blow-up of a map tangent to the identity, and
$[v]\in E$ is a non-degenerate characteristic direction with non-zero
director~$\alpha$, then $\iota_{[v]}(f,E)=1/\alpha$. 

\Rem Theorem~\index\ is only a very particular case of a much more
general index theorem, valid for holomorphic self-maps of complex manifolds of
any dimension fixing pointwise a smooth complex submanifold of any codimension,
or a hypersurface even with singularities; see [BrT], [Br] and [ABT], where some
applications to dynamics are also discussed. In particular, in [ABT] is
introduced a canonical section of a suitable vector bundle describing the local
dynamics in an infinitesimal neighbourhood of the submanifold, providing
in particular a more intrinsic description of the index.

Now, a combinatorial argument (again inspired by Camacho and Sad~[CS]) shows
that if we have $f\in\End(\C^2,O)$ with an isolated fixed point, and such that
applying the blow-up procedure to the lifted map~$\tilde f$ starting from a
singular direction~$[v]\in\P^1(\C)=E$ we end up with only $(\star_2)$
singularities, then the index of~$\tilde f$ at~$[v]$ along~$E$ must be a
non-negative rational number. But the first Chern class of~$N_E$ is~$-1$, and so
there must be at least one singular directions whose index is not a non-negative
rational number, and thus the blow-up procedure must yield at least one
dicritical or $(\star_1)$ singularity, and hence a Fatou flower for our map~$f$,
completing the proof of Theorem~\Abate.

Actually, we have proved the following slightly more precise result:

\newthm Theorem \Abatedue: (Abate, 2001 [A2]) Let $f\in\End(\C^2,O)$ be a
holomorphic local dynamical system tangent to the identity and with an isolated
fixed point at the origin. Let $[v]\in\P^1(\C)$ be a singular direction such
that~$\iota_{[v]}\bigl(\tilde f,\P^1(\C)\bigr)\notin\Q^+$, where $\tilde f$ is
the lift of~$f$ to the blow-up of the origin. Then $f$ has a Fatou flower
tangent to~$[v]$.

\Rem To be even more precise, Theorem~\Abatedue\ is more a statement on the
lifted map~$\tilde f$ than on the original~$f$. Indeed, the argument used to
prove Theorem~\Abatedue\ (or a similar argument along the lines of~[Ca]) can be
used to prove the following: {\it let $E$ be a (not necessarily compact) Riemann
surface inside a 2-dimensional complex manifold~$M$, and take $f\in\End(M,E)$
tangent to the identity at all points of~$E$ and tangential to~$E$. Let $p\in E$
be a singular point of~$f$ such that $\iota_p(f,E)\notin\Q^+$. Then there exist
parabolic curves for~$f$ at~$p$.} This latter statement has been recently
generalized in two ways. Degli Innocenti [DI] has proved that we can allow~$E$
to be singular at~$p$ (but irreducible; in the reducible case one has to impose
conditions on the indeces of~$f$ along all irreducible components of~$E$ passing
through~$p$).  Molino [Mo], on the other hand, has proved that the statement
still holds assuming only $\iota_p(f,E)\ne 0$, at least for $f$ of order~2 (and
$E$ smooth at~$p$); it is natural to conjecture that this should be true for $f$
of any order.

As already remarked, the reduction of singularities via
blow-ups seem to work only in dimension~2. This leaves open the problem of the
validity of something like Theorem~\Abate\ in dimension~$n\ge 3$; see [AT]
for some partial results.

Furthermore, as far as I know, it is completely open, even in dimension~2, the
problem of describing the stable set of a holomorphic local dynamical
system tangent to the identity, as well as the more general problem of the
topological classification of such dynamical systems. Some results in
the case of a dicritical singularity are presented in~[BM]. 

We end this section with a couple of words on holomorphic local dynamical
systems with a parabolic fixed point where the differential is not
diagonalizable. Particular examples are studied in detail in~[CD], [A4]
and~[GS]. In~[A1] it is described
a canonical procedure for lifting an~$f\in\End(\C^n,O)$ whose differential at the
origin is not diagonalizable to a map defined in a suitable iterated blow-up of
the origin (obtained blowing-up not only points but more general submanifolds)
with a canonical fixed point where the differential is diagonalizable. Using
this procedure it is for instance possible to prove the following

\newthm Corollary \Jordan: (Abate, 2001 [A2]) Let $f\in\End(\C^2,O)$ be a
holomorphic local dynamical system with $df_O=J_2$, the canonical Jordan matrix
associated to the eigenvalue~$1$, and assume that the origin is an isolated
fixed point. Then $f$ admits at least one parabolic curve tangent to~$[1:0]$ at
the origin. 

\smallsect 7. Several complex variables: other cases

Outside the hyperbolic and parabolic cases, there are not that many general
results. Theorems~\Bryunohyp\ and~\BryunoCremer\ apply to the elliptic case too,
but, as already remarked, it is not known whether the Bryuno condition is still
necessary for holomorphic linearizability, that is, if any analogue of
Theorem~\BY.(ii) holds in several variables. However,
another result in the spirit of Theorem~\BryunoCremer\ is the following:

\newthm Theorem \YC: (Yoccoz, 1995 [Y2]) Let $A\in GL(n,\C)$ be an invertible
matrix such that its eigenvalues have no resonances and such that its Jordan
normal form contains a non-trivial block associated to an eigenvalue of modulus
one. Then there exists $f\in\End(\C^n,O)$ with $df_O=A$ which is not
holomorphically linearizable.

A case that has received some attention is the so-called
semi-attractive case: a holomorphic local dynamical system $f\in\End(\C^n,O)$
is said {\sl semi-attractive} if the eigenvalues of~$df_O$ are either equal to~1
or with modulus strictly less than~1. The dynamics of semi-attractive dynamical
systems has been studied in detail by Fatou~[F4], Nishimura~[N],
Ueda~[U1--2], Hakim~[H1] and Rivi~[Ri--2]. Their results more or less
say that the eigenvalue~1 yields the existence of a ``parabolic manifold"~$M$ ---
in the sense of Theorem~\Hakim.(ii) --- of a suitable dimension, while the
eigenvalues with modulus less than one ensure, roughly speaking, that the
orbits of points in the normal bundle of~$M$ close enough to~$M$ are attracted
to it. For instance, Rivi proved the following

\newthm Theorem \Rivi: (Rivi, 1999 [Ri1--2]) Let $f\in\End(\C^n,O)$ be a
holomorphic local dynamical system. Assume that $1$ is an eigenvalue of
(algebraic and geometric) multiplicity~$q\ge 1$ of $df_O$, and that the other
eigenvalues of~$df_O$ have modulus less than~$1$. Then:
{\smallskip
\item{\rm(i)}We can choose local coordinates $(z,w)\in\C^q\times\C^{n-q}$ such
that $f$ expressed in these coordinates becomes
$$
\cases{f_1(z,w)=A(w)z+P_{2,w}(z)+P_{3,w}(z)+\cdots,\cr
\noalign{\smallskip}
f_2(z,w)=G(w)+B(z,w)z,\cr}
$$
where: $A(w)$ is a $q\times q$ matrix with entries holomorphic in~$w$ and
$A(O)=I_q$; the $P_{j,w}$ are $q$-uples of homogeneous polynomials in~$z$ of
degree~$j$ whose coefficients are holomorphic in~$w$; $G$ is a holomorphic
self-map of~$\C^{n-q}$ such that $G(O)=O$ and the eigenvalues of~$dG_O$ are
the eigenvalues of~$df_O$ with modulus strictly less than~$1$; and $B(z,w)$ is
an $(n-q)\times q$ matrix of holomorphic functions vanishing at the origin. In
particular, $f_1(z,O)$ is tangent to the identity.
\item{\rm(ii)} If $v\in\C^q\subset\C^m$ is a non-degenerate characteristic
direction for~$f_1(z,O)$, and the latter map has order~$\nu$,\break\indent then
there exist
$\nu-1$ disjoint $f$-invariant $(n-q+1)$-dimensional complex
submanifolds~$M_j$ of~$\C^n$, with\break\indent the origin in their boundary,
such that the orbit of every point of~$M_j$ converges to the origin tangentially
\break\indent to~$\C v\oplus E$, where $E\subset\C^n$ is the subspace generated
by the generalized eigenspaces associated to the\break\indent eigenvalues
of~$df_O$ with modulus less than one.}  

Rivi also has results in the spirit of Theorem~\Hakim, and results when the
algebraic and geometric multiplicities of the eigenvalue~$1$ differ; see~[Ri1, 2]
for details.

As far as I know, there are no results on the formal or holomorphic
classification of semi-attractive holomorphic local dynamical systems. However,
Canille Martins has given the topological classification in dimension~2, using
Theorem~\Camacho\ and general results on normally hyperbolic dynamical systems
due to Palis and Takens~[PT]:

\newthm Theorem \CanilleMartins: (Canille Martins, 1992 [CM]) Let
$f\in\End(\C^2,O)$ be a holomorphic local dynamical system such that $df_O$ has
two eigenvalues $\lambda_1$,~$\lambda_2\in\C$, where $\lambda_1$ is a
primitive $q$-root of unity, and $|\lambda_2|\ne 0,1$. Then $f$ is topologically
locally conjugated to the map
$$
(z,w)\mapsto (\lambda_1 z+z^{kq+1},\lambda_2 w)
$$
for a suitable $k\in\N^*$.

We end this survey by recalling a very recent result by Bracci and Molino. Assume
that $f\in\End(\C^2,O)$ is a holomorphic local dynamical system such that the
eigenvalues of~$df_O$ are~$1$ and~$e^{2\pi i\theta}\ne 1$. If $e^{2\pi i\theta}$
satisfies the Bryuno condition, P\"oschel~[P\"o] proved the existence of a
1-dimensional $f$-invariant holomorphic disk containing the origin where~$f$ is
conjugated to the irrational rotation of angle~$\theta$. On the other hand,
Bracci and Molino give sufficient conditions (depending on~$f$ but not
on~$e^{2\pi i\theta}$, expressed in terms of two new holomorphic invariants, and
satisfied by generic maps) for the existence of parabolic curves tangent to the
eigenspace of the eigenvalue~1; see [BrM] for details.

\setref{EHRS}
\beginsection References

\art A1 M. Abate: Diagonalization of non-diagonalizable discrete holomorphic
dynamical systems! Amer. J. Math.! 122 2000 757-781

\art A2 M. Abate: The residual index and the dynamics of holomorphic maps
tangent to the identity! Duke Math. J.! 107 2001 173-207

\book A3 M. Abate: An introduction to hyperbolic dynamical systems! I.E.P.I.
Pisa, 2001

\art A4 M. Abate: Basins of attraction in quadratic dynamical systems with a
Jordan fixed point! Nonlinear Anal.! 51 2002 271-282

\art AT M. Abate, F. Tovena: Parabolic curves in $\C^3$! 
Abstr. Appl. Anal.! 2003 2003 275-294

\pre ABT M. Abate, F. Bracci, F. Tovena: Index theorems for holomorphic
self-maps! To appear in Ann. of Math.! 2003

\book Ar V.I. Arnold: Geometrical methods in the theory of ordinary differential
equations! Springer-Verlag, Berlin, 1988

\art B L.E. B\"ottcher: The principal laws of convergence of iterates and their 
application to analysis! Izv. Kazan. Fiz.-Mat. Obshch.! 14 1904
155-234

\pre Br F. Bracci: The dynamics of holomorphic maps near curves of fixed points!
To appear in Ann. Scuola Norm. Sup. Pisa! 2002

\pre BrM F. Bracci, L. Molino: The dynamics near quasi-parabolic fixed points of 
holomorphic diffeomorphisms in $\C^2$! To appear in Amer. J. Math.! 2003

\pre BrS F. Bracci, T. Suwa: Residues for singular pairs and dynamics of
biholomorphic maps of singular surfaces! Preprint! 2003

\art BrT F. Bracci, F. Tovena: Residual indices of
holomorphic maps relative to  singular curves of fixed points on
surfaces! Math. Z.! 242 2002 481-490

\art BM F.E. Brochero Mart\'\i nez: Groups of germs of analytic diffeomorphisms
in $(\C^2,O)$! J. Dynamic. Control Systems! 9 2003 1-32

\art Bry1 A.D. Bryuno: Convergence of transformations of differential equations
to normal forms! Dokl. Akad. Nauk. USSR! 165 1965 987-989

\art Bry2 A.D. Bryuno: Analytical form of differential equations, I! Trans.
Moscow  Math. Soc.! 25 1971 131-288

\art Bry3 A.D. Bryuno: Analytical form of differential equations, I\negthinspace
I! Trans. Moscow  Math. Soc.! 26 1972 199-239

\art C C. Camacho: On the local structure of conformal mappings and holomorphic
vector fields! Ast\'e\-risque! 59--60 1978 83-94

\art CS C. Camacho, P. Sad: Invariant varieties through singularities of holomorphic vector
fields! Ann. of Math.! 115 1982 579-595

\art CM J.C. Canille Martins: Holomorphic flows in $(\C^3,O)$ with resonances!
Trans. Am. Math. Soc.! 329 1992 825-837

\art Ca J. Cano: Construction of invariant curves for singular holomorphic
vector fields! Proc. Am. Math. Soc.! 125 1997 2649-2650

\book Ch M. Chaperon: G\'eom\'etrie diff\'erentielle et singularit\'es des
syst\`emes dynamiques! Ast\'e\-risque {\bf 138--139,} 1986

\art CD D. Coman, M. Dabija: On the dynamics of some diffeomorphisms of $\C^2$
near parabolic fixed points! Houston J. Math.! 24 1998 85-96

\art Cr1 H. Cremer: Zum Zentrumproblem! Math. An..! 98 1927 151-163

\art Cr2 H. Cremer: \"Uber die H\"aufigkeit der Nichtzentren! Math. Ann.! 115
1938 573-580

\book DI F. Degli Innocenti: Dinamica di germi di foliazioni e diffeomorfismi
olomorfi vicino a curve singolari! Undergraduate thesis, Universit\`a di Firenze,
2003

\art DG D. DeLatte, T. Gramchev: Biholomorphic maps with linear parts having
Jordan blocks: linearization and resonance type phenomena! Math. Phys. El. J.! 8
2002 1-27

\book \'E1 J. \'Ecalle: Les fonctions r\'esurgentes. Tome
I: Les alg\`ebres de fonctions r\'esurgentes! Publ. Math. Orsay {\bf 81-05,}
Universit\'e de Paris-Sud, Orsay, 1981 

\book \'E2 J. \'Ecalle: Les fonctions r\'esurgentes. Tome
I\negthinspace I: Les fonctions r\'esurgentes appliqu\'ees \`a l'it\'eration!
Publ. Math. Orsay {\bf 81-06,} Universit\'e de Paris-Sud, Orsay, 1981 

\book \'E3 J. \'Ecalle: Les fonctions r\'esurgentes. Tome I\negthinspace
I\negthinspace I: L'\'equation du pont et la classification analytique des
objects locaux! Publ. Math. Orsay {\bf 85-05,} Universit\'e de Paris-Sud, Orsay,
1985 

\coll \'E4 J. \'Ecalle: Iteration and analytic classification of local
diffeomorphisms of $\C^\nu$! Iteration theory and its functional equations!
Lect. Notes in Math. {\bf 1163,} Springer-Verlag, Berlin, 1985, pp. 41--48

\art \'EV J. \'Ecalle, B. Vallet: Correction and linearization of resonant 
vector fields and diffeomorphisms! Math. Z.! 229 1998 249-318

\art F1 P. Fatou: Sur les \'equations fonctionnelles, I! Bull. Soc. Math. France!
47 1919 161-271

\art F2 P. Fatou: Sur les \'equations fonctionnelles, I\negthinspace I! Bull.
Soc. Math. France! 48 1920 33-94

\art F3 P. Fatou: Sur les \'equations fonctionnelles, I\negthinspace
I\negthinspace I! Bull. Soc. Math. France! 48 1920 208-314

\art F4 P. Fatou: Substitutions analytiques et \'equations fonctionnelles \`a
deux variables! Ann. Sc. Ec. Norm. Sup.! 40 1924 67-142

\pre FJ C. Favre, M. Jonsson: The valuative tree! Preprint, 
arXiv: math.AC/0210265! 2002

\coll FHY A. Fathi, M. Herman, J.-C. Yoccoz: A proof of Pesin's stable manifold
theorem! Geometric Dynamics! Lect Notes in Math. 1007, Springer Verlag, Berlin,
1983, pp.~177-216

\art GS V. Gelfreich, D. Sauzin: Borel summation and splitting of separatrices
for the H\'enon map! Ann. Inst. Fourier Grenoble! 51 2001 513-567

\art G1 D.M. Grobman: Homeomorphism of systems of differential equations! Dokl.
Akad. Nauk. USSR! 128 1959 880-881

\art G2 D.M. Grobman: Topological classification of neighbourhoods of a
singularity in $n$-space! Math. Sbornik! 56 1962 77-94

\art H J.S. Hadamard: Sur l'it\'eration et les solutions asymptotyques des
\'equations diff\'erentielles! Bull. Soc. Math. France! 29 1901 224-228

\art Ha1 M. Hakim: Attracting domains for semi-attractive transformations
of~$\C^p$! Publ. Matem.! 38 1994 479-499

\art Ha2 M. Hakim: Analytic transformations of $(\C^p,0)$ tangent to the
identity! Duke Math. J.! 92 1998 403-428

\pre Ha3 M. Hakim: Transformations tangent to the identity. Stable pieces
of manifolds! Preprint! 1997

\art Har P. Hartman: A lemma in the theory of structural stability of
differential equations! Proc. Am. Math. Soc.! 11 1960 610-620

\book HK B. Hasselblatt, A. Katok: Introduction to the modern theory of
dynamical systems! Cambridge Univ. Press, Cambridge, 1995

\coll He M. Herman: Recent results and some open questions on Siegel's
linearization theorem of germs of complex analytic diffeomorphisms of $\C^n$
near a fixed point! Proc. $8^{th}$ Int. Cong. Math. Phys.! World Scientific,
Singapore, 1986, pp. 138--198

\book HPS M. Hirsch, C.C. Pugh, M. Shub: Invariant manifolds! Lect. Notes Math.
{\bf 583,} Springer-Verlag, Berlin, 1977

\art HP J.H. Hubbard, P. Papadopol: Superattractive fixed points in $\C^n$!
Indiana Univ. Math. J.! 43 1994 321-365

\art I1 Yu.S. Il'yashenko: Divergence of series reducing an analytic differential
equation to linear normal form at a singular point! Funct. Anal. Appl.! 13 1979
227-229

\coll I2 Yu.S. Il'yashenko: Nonlinear Stokes phenomena! Nonlinear Stokes
phenomena! Adv. in Soviet Math. {\bf 14,} Am. Math. Soc., Providence, 1993,
pp.~1--55

\art K T. Kimura: On the iteration of analytic functions! Funk. Eqvacioj! 14
1971 197-238

\art {K\oe}  G. K\oe nigs: Recherches sur les integrals de certain equations
fonctionelles! Ann. Sci. \'Ec. Norm. Sup.! 1 1884 1-41

\art L L. Leau: \'Etude sur les equations fonctionelles \`a une ou plusieurs variables! Ann.
Fac. Sci. Toulouse! 11 1897 E1-E110

\art M1 B. Malgrange: Travaux d'\'Ecalle et de Martinet-Ramis sur les syst\`emes 
dynamiques! Ast\'erisque! 92-93 1981/82 59-73

\art M2 B. Malgrange: Introduction aux travaux de J. \'Ecalle! Ens. Math.!
31 1985 261-282

\art MM J.F. Mattei, R. Moussu: Holonomie et int\'egrales premi\`eres! Ann. Scient. Ec.
Norm. Sup.! 13 1980 469-523

\book Mi J. Milnor: Dynamics in one complex variable! Vieweg, Braunschweig, 2000

\pre Mo L. Molino: Parabolic curves at singular points! In preparation! 2003

\art N Y. Nishimura: Automorphismes analytiques admettant des sousvari\'et\'es
de point fixes attractives dans la direction transversale! J. Math. Kyoto Univ.!
23 1983 289-299

\art PT J. Palis, F. Takens: Topological equivalence of normally hyperbolic
dynamical systems! Topology! 16 1977 335-345

\art P1 R. P\'erez-Marco: Sur les dynamiques holomorphes non lin\'earisables et
une conjecture de V.I. Arnold! Ann. Sci. \'Ecole Norm. Sup.! 26 1993 565-644

\pre P2 R. P\'erez-Marco: Topology of Julia sets and hedgehogs! Preprint!
Universit\'e de Paris-Sud, 94-48, 1994

\art P3 R. P\'erez-Marco: Non-linearizable holomorphic dynamics having an
uncountable number of symmetries! Invent. Math.! 199 1995 67-127

\pre P4 R. P\'erez-Marco: Holomorphic germs of quadratic type! Preprint! 1995

\pre P5 R. P\'erez-Marco: Hedgehogs dynamics! Preprint! 1995

\art P6 R. P\'erez-Marco: Sur une question de Dulac et Fatou! C.R. Acad. Sci.
Paris! 321 1995 1045-1048

\art P7 R. P\'erez-Marco: Fixed points and circle maps! Acta Math.! 179 1997
243-294

\pre P8 R. P\'erez-Marco: Linearization of holomorphic germs with resonant
linear part! Preprint, arXiv: math.DS/0009030! 2000

\pre P9 R. P\'erez-Marco: Total convergence or general divergence in small
divisors! Preprint, arXiv: math. DS/0009029! 2000

\art Pe O. Perron: \"Uber Stabilit\"at und asymptotisches Verhalten der Integrale
von Differentialgleichungssystemen! Math. Z.! 29 1928 129-160

\art Pes Ja.B. Pesin: Families of invariant manifolds corresponding to non-zero
characteristic exponents! Math. USSR Izv.! 10 1976 1261-1305

\book Po H. Poincar\'e: \OE uvres, Tome I! Gauthier-Villars, Paris, 1928,
pp.~XXXVI--CXXIX

\art P\"o J. P\"oschel: On invariant manifolds of complex analytic mappings near
fixed points! Exp. Math.! 4 1986 97-109

\art R1 L. Reich: Das Typenproblem bei formal-biholomorphien Abbildungen mit
anziehendem Fixpunkt! Math. Ann.! 179 1969 227-250

\art R2 L. Reich: Normalformen biholomorpher Abbildungen mit anziehendem
Fixpunkt! Math. Ann.! 180 1969 233-255

\book Ri1 M. Rivi: Local behaviour of discrete dynamical systems! Ph.D. Thesis,
Universit\`a di Firenze, 1999

\art Ri2 M. Rivi: Parabolic manifolds for semi-attractive holomorphic germs!
Mich. Math. J.! 49 2001 211-241

\art S A.A. Shcherbakov: Topological classification of germs of conformal
mappings with identity linear part! Moscow Univ. Math. Bull.! 37 1982 60-65

\book Sh M. Shub: Global stability of dynamical systems! Springer, Berlin, 1987

\art Si C.L. Siegel: Iteration of analytic functions! Ann. of Math.! 43 1942
607-612

\art St1 S. Sternberg: Local contractions and a theorem of Poincar\'e! Amer. J.
Math.! 79 1957 809-824

\art St2 S. Sternberg: The structure of local homomorphisms, I\negthinspace I!
Amer. J. Math.! 80 1958 623-631

\art U1 T. Ueda: Local structure of analytic transformations of two complex
variables, I! J. Math. Kyoto Univ.! 26 1986 233-261

\art U2 T. Ueda: Local structure of analytic transformations of two complex
variables, I\negthinspace I! J. Math. Kyoto Univ.! 31 1991 695-711

\art Us S. Ushiki: Parabolic fixed points of two-dimensional complex dynamical
systems! S\=urikaisekiken\-ky\=usho K\=oky\=uroku! 959 1996 168-180

\art V S.M. Voronin: Analytic classification of germs of conformal maps
$(\C,0)\to (\C,0)$ with identity linear part! Func. Anal. Appl.! 15 1981 1-17

\art Y1  J.-C. Yoccoz: Lin\'earisation des germes de diff\'eomorphismes
holomorphes de $(\C,0)$! C.R. Acad. Sci. Paris! 306 1988 55-58

\art Y2 J.-C. Yoccoz: Th\'eor\`eme de Siegel, nombres de Bryuno et polyn\^omes
quadratiques! Ast\'e\-ris\-que! 231 1995 3-88

\art W B.J. Weickert: Attracting basins for automorphisms of $\C^2$! Invent.
Math.! 132 1998 581-605

\art Wu H. Wu: Complex stable manifolds of holomorphic diffeomorphisms! Indiana
Univ. Math. J.! 42 1993 1349-1358

\bye